\documentclass{CSML}

\def\dOi{11(4:6)2015}
\lmcsheading%
{\dOi}
{1--30}
{}
{}
{Nov.~11, 2013}
{Dec.~\phantom02, 2015}
{}

\ACMCCS{[{\bf Theory of computation}]: Logic; [{\bf Mathematics of
    computing}]:  Mathematical analysis---Nonlinear equations}

\usepackage{hyperref}

\usepackage{amssymb}
\usepackage{tikz}
\usetikzlibrary{arrows,chains,matrix,positioning,scopes}
\makeatletter
\tikzset{join/.code=\tikzset{after node path={%
\ifx\tikzchainprevious\pgfutil@empty\else(\tikzchainprevious)%
edge[every join]#1(\tikzchaincurrent)\fi}}}
\makeatother
\tikzset{>=stealth',every on chain/.append style={join},
         every join/.style={->}}
\tikzstyle{labeled}=[execute at begin node=$\scriptstyle,
   execute at end node=$]

\theoremstyle{definition}
\newtheorem{theorem}{Theorem}

\newtheorem{observation}[theorem]{Observation}
\newtheorem{question}[thm]{Open Question}

\newcommand{\dom}{\operatorname{dom}}
\newcommand{\id}{\textnormal{id}}
\newcommand{\Cantor}{{\{0, 1\}^\mathbb{N}}}
\newcommand{\Baire}{{\mathbb{N}^\mathbb{N}}}

\newcommand{\hide}[1]{}
\newcommand{\mto}{\rightrightarrows}
\newcommand{\uint}{{[0, 1]}}

\newcommand{\name}[1]{\textsc{#1}}
\newcommand{\C}{\textrm{C}}
\newcommand{\ic}[1]{\textrm{C}_{\sharp #1}}
\newcommand{\xc}[1]{\textrm{XC}_{#1}}

\newcommand{\etal}{et al.~}
\newcommand{\Sierp}{Sierpi\'nski }
\newcommand{\leqW}{\leq_{\textrm{W}}}
\newcommand{\leW}{<_{\textrm{W}}}
\newcommand{\equivW}{\equiv_{\textrm{W}}}

\newcommand{\nleqW}{\nleq_\textrm{W}}

\begin{document}

\title{Finite choice, convex choice and finding roots}

\author[S.~Le Roux]{St\'ephane Le Roux\rsuper a}	%required
\address{{\lsuper a}Computer Science departement\\ Universit\'e Libre de Bruxelles, Belgium}	%required
\email{Stephane.Le.Roux@ulb.ac.be}  %optional
\thanks{Le Roux was at Technische Universit\"at Darmstadt, department of Mathematics, when this research started.}	%optional

\author[A.~Pauly]{Arno Pauly\rsuper b}	%optional
\address{{\lsuper b}Clare College\\ University of Cambridge, United Kingdom}	%optional
\email{Arno.Pauly@cl.cam.ac.uk}  %optional

%% required for running head on odd and even pages, use suitable
%% abbreviations in case of long titles and many authors:

%% mandatory lists of keywords and classifications:
\titlecomment{{\lsuper*}An extended abstract appeared in the proceedings of CiE 2013 \cite{paulyleroux-cie}.}
%%%%%%%%%%%%%%%%%%%%%%%%%%%%%%%%%%%%%%%%%%%%%%%%%%%%%%%%%%%%%%%%%%%%%%%%%%%

%% the abstract has to PRECEED the command \maketitle:
%% be sure not to issue the \maketitle command twice!

\begin{abstract}
  \noindent We investigate choice principles in the Weihrauch lattice for finite sets on the one hand, and convex sets on the other hand. Increasing cardinality and increasing dimension both correspond to increasing Weihrauch degrees. Moreover, we demonstrate that the dimension of convex sets can be characterized by the cardinality of finite sets encodable into them. Precisely, choice from an $n + 1$ point set is reducible to choice from a convex set of dimension $n$, but not reducible to choice from a convex set of dimension $n - 1$.

Furthermore we consider searching for zeros of continuous functions. We provide an algorithm producing $3^n$ real numbers containing all zeros of a continuous function with up to $n$ local minima. This demonstrates that having finitely many zeros is a strictly weaker condition than having finitely many local extrema. We can prove $3^n$ to be optimal.
\end{abstract}

\maketitle

\section{Introduction}
In the investigation of the computational content of mathematical theorems in the Weihrauch lattice, variations of closed choice principles have emerged as useful canonical characterizations \cite{brattka3, paulybrattka, paulybrattka3cie}. Closed choice principles are multivalued functions taking as input a non-empty closed subset of some fixed space, and have to provide some element of the closed set as output. In \cite{brattka3, paulybrattka} the influence of the space on the computational difficulty of (full) closed choice was investigated, whereas in \cite{paulybrattka3cie} it turned out that the restriction of choice to connected closed subsets of the unit hypercube is equivalent to Brouwer's Fixed Point theorem for the same space.

Here the restrictions of closed choice to convex subsets (of the unit hypercube of dimension $n$), and to finite subsets (of a compact metric space) are the foci of our investigations. Via the connection between closed choice and non-deterministic computation \cite{ziegler2, paulybrattka, paulybrattka2, paulydebrecht}, in particular the latter problem is prototypic for those problems having only finitely many correct solutions where wrong solutions are identifiable. As such, some parts may be reminiscent of some ideas from \cite{kreinovich, gasarch}.

One of our main results shows that choice for finite sets of cardinality $n + 1$ can be reduced to choice for convex sets of dimension $n$, but not to convex choice of dimension $n-1$. This demonstrates a computational aspect in which convex sets get more complicated with increasing dimension. As such, our work also continues the study of the structural complexity of various classes of subsets of the unit hypercubes done in \cite{leroux, kihara}.

Some of the techniques used to establish our main results are promising with regards to further applicability to other classes of choice principles, or to even more general Weihrauch degrees. These techniques are presented in Section \ref{sec:relativetechniques}.

Finally, some of the results are transferred to the problem of finding zeros of a continuous function. We show that finding zeros of a function merely guaranteed to have finitely many zeros is strictly harder than finding zeros of a function having finitely many zeros as a consequence of having finitely many local extrema. This is achieved via an algorithm that lists $3^n$ potential zeros of a continuous function $f : \uint \to \mathbb{R}$ with up to $n$ local minima and $f(0) > 0 \wedge f(1) \neq 0$ while guaranteeing that all true zeros are listed. We show that any algorithm of this kind needs at least $3^n$ guesses to cover all zeros.

\subsection{Weihrauch reducibility}
We briefly recall some basic results and definitions regarding
the Weihrauch lattice. The original definition of Weihrauch reducibility is due to Weihrauch
and has been studied for many years (see \cite{stein,weihrauchb,weihrauchc,mylatz,hertling,paulymaster,mylatzb}).
Rather recently it has been noticed that a certain variant of this reducibility yields
a lattice that is very suitable for the classification of mathematical theorems
(see \cite{gherardi,brattka2,brattka3,paulyreducibilitylattice,paulybrattka,paulyincomputabilitynashequilibria,gherardi4,paulykojiro,shafer,paulyleroux3-arxiv,pauly-fouche}). A basic reference for notions
from computable analysis is \cite{weihrauchd}.
The Weihrauch lattice is a lattice of multi-valued functions on represented
spaces. A represented space is a pair $(X, \delta_X)$ where $\delta_X : \subseteq \Baire \to X$ is a partial surjection, called representation.
In general we use the symbol ``$\subseteq$'' in order to indicate that a function is potentially partial.
Using represented spaces we can define the concept of a realizer. We denote the composition of
two (multi-valued) functions $f$ and $g$ either by $f\circ g$ or by $fg$.

\begin{defi}[Realizer]
Let $f : \subseteq (X, \delta_X) \mto (Y, \delta_Y)$ be a multi-valued function on represented spaces.
A function $F:\subseteq\Baire\to\Baire$ is called a {\em realizer} of $f$, in symbols $F\vdash f$, if
$\delta_YF(p)\in f\delta_X(p)$ for all $p\in\dom(f\delta_X)$.
\end{defi}

Realizers allow us to transfer the notions of computability, continuity, and other notions available for Baire space to any represented space;
a function between represented spaces will be called {\em computable}, if it has a computable realizer, etc.

\label{page:fa}
We will need a generalization of the restriction of a multivalued function to a subspace of its domain. For some represented space $\mathbf{X} = (X, \delta_X)$ and $A \subseteq \Baire$, we use $\mathbf{X}_A$ to denote the represented space $(\delta_X[A], (\delta_X)|_{A})$. This is a proper generalization of the notion of a subspace. Given $f :\subseteq \mathbf{X} \mto \mathbf{Y}$ and $A \subseteq \Baire$, then $f_A$ is the induced map $f_A :\subseteq \mathbf{X}_A \mto \mathbf{Y}$.

Now we can define Weihrauch reducibility, using $\langle,\rangle$ to denote some standard pairing on Baire space.

\begin{defi}[Weihrauch reducibility]
Let $f,g$ be multi-valued functions on represented spaces.
Then $f$ is said to be {\em Weihrauch reducible} to $g$, in symbols $f\leqW g$, if there are computable
functions $K,H:\subseteq\Baire\to\Baire$ such that $K\langle \id, GH \rangle \vdash f$ for all $G \vdash g$.
Moreover, $f$ is said to be {\em strongly Weihrauch reducible} to $g$, in symbols $f\leq_{sW} g$,
if there are computable functions $K,H$ such that $KGH\vdash f$ for all $G\vdash g$.
\end{defi}

We note that the relations $\leqW$, $\leq_{sW}$ and $\vdash$ implicitly refer to the underlying representations, which
we mention explicitly only when necessary. It is known that these relations only depend on the underlying equivalence
classes of representations, but not on the specific representatives (see Lemma~2.11 in \cite{brattka2}).
We use $\equivW$ and $\equiv_{sW}$ to denote the respective equivalences regarding $\leqW$ and $\leq_{sW}$,
and by $\leW$ and $<_{sW}$ we denote strict reducibility.

There are three operations defined on Weihrauch degrees that are used in the present paper, product $\times$, composition $\star$ and Kleene star $^*$. The former operation was originally introduced in \cite{brattka2,paulyreducibilitylattice}, the second in \cite{gherardi4} and the third in \cite{paulyreducibilitylattice,paulyincomputabilitynashequilibria}. Informally, access to the product allows us to use both involved operations independently, whereas for the composition the call to the left operation may depend on the answer received from the right. The Kleene star corresponds to any finite number of parallel uses.

\begin{defi}
Given $f :\subseteq \mathbf{X} \mto \mathbf{Y}$, $g : \subseteq \mathbf{U} \mto \mathbf{V}$, define $f \times g : \subseteq (\mathbf{X} \times \mathbf{U}) \mto (\mathbf{Y} \times \mathbf{V})$ via $(y, v) \in (f \times g)(x, u)$ iff $y \in f(x)$ and $v \in g(u)$.
\end{defi}

\begin{defi}
Given $f :\subseteq \mathbf{X} \mto \mathbf{Y}$, $g : \subseteq \mathbf{U} \mto \mathbf{V}$, let \[f \star g := \sup_{\leqW} \{f' \circ g' \mid f' \leqW f \wedge g' \leqW g\}\]
where $f'$, $g'$ are understood to range over all those multivalued functions where the composition is defined.
\end{defi}

Both $\times$ and $\star$ are associative, but only $\times$ is commutative. We point out that while it is not obvious that the supremum in the definition of $\star$ always exists, this is indeed the case, hence $\star$ is actually a total operation. We will iterate both $\times$ and $\star$, writing $f^0 = f^{(0)} = \id_\Baire$ and $f^{n+1} = f \times f^{n}$, $f^{(n+1)} = f \star f^{(n)}$. The former is subsequently used to introduce $f^*$ via $f^*(n,(x_1, \ldots, x_n)) = f^n(x_1, \ldots, x_n)$.

We will also refer to a special Weihrauch degree, denoted by $0$. Its representatives are the nowhere defined functions, and it is the bottom element of the lattice.
\subsection{Closed Choice and variations thereof}
The space of continuous functions from a represented space $\mathbf{X}$ to $\mathbf{Y}$ has a natural representation itself, as a consequence of the \textrm{UTM}-theorem. This represented space is denoted by $\mathcal{C}(\mathbf{X}, \mathbf{Y})$.

A special represented space of utmost importance is \Sierp space $\mathbb{S}$ containing two elements $\{\top, \bot\}$ represented by $\delta_\mathbb{S} : \Baire \to \mathbb{S}$ where $\delta_\mathbb{S}(0^\mathbb{N}) = \bot$ and $\delta_\mathbb{S}(p) = \top$, iff $p \neq 0^\mathbb{N}$. The space $\mathcal{A}(\mathbf{X})$ of closed subsets of $\mathbf{X}$ is obtained from $\mathcal{C}(\mathbf{X}, \mathbb{S})$ by identifying a set $A \subseteq \mathbf{X}$ with the characteristic function $\chi_{X \setminus A} : \mathbf{X} \to \mathbb{S}$ of its complement.

For a computable metric space $\mathbf{X}$, an equivalent representation $\psi_-:\Baire\to \mathcal{A}(\mathbf{X})$, can be defined by
$\psi_-(p):=X\setminus\bigcup_{i=0}^\infty B_{p(i)}$,
where $(B_n)_{n \in \mathbb{N}}$ is some standard enumeration of the open balls of $X$ with centers in the dense subset and rational radii (possibly $0$).
The computable points in $\mathcal{A}(X)$ are called {\em co-c.e.\ closed sets}.
We are primarily interested in closed subsets of computable metric spaces; additionally, most of our considerations pertain to compact spaces -- see Subsection \ref{subsec:beyondcompactness} for the exceptions.

The computability structure available on the closed sets mostly follows the intuitive expectations, for an explicit treatment we refer primarily to \cite{pauly-synthetic-arxiv}. We will also need computability of the closed convex hull, which we shall use for subspaces of $\mathbb{R}^n$, in short, Euclidean spaces. Recall from \cite{pauly-synthetic-arxiv} that a space is called \emph{computably compact}, iff it is semidecidable (recognizable) that an open set contains a closed set.
\begin{prop}[\footnote{This has essentially already been observed by \name{Ziegler} \cite{ziegler8}.}]
\label{prop:convexhull}
Let $\mathbf{X}$ be a computably compact\footnote{Compactness is a
  necessary condition here: $\overline{\operatorname{ConvexHull}} :
  \mathcal{A}(\mathbb{R}) \to \mathcal{A}(\mathbb{R})$ is not
  computable. \begin{proof}Assume the contrary. Given some $A \in
    \mathcal{A}(\mathbb{R})$, compute $A_1 :=
    \overline{\operatorname{ConvexHull}}(\{-1\} \cup ([0;\infty) \cap
    A))$ and $A_2 := \overline{\operatorname{ConvexHull}}(\{1\} \cup
    ((-\infty;0] \cap A))$. Now note
    $\operatorname{IsEmpty}_\mathbb{R}(A) = (-\frac{1}{2} \notin A_1)
    \wedge (\frac{1}{2} \notin A_2)$ -- but computability of
    $\operatorname{IsEmpty}_\mathbb{R}$ would imply compactness of
    $\mathbb{R}$.\end{proof}} Euclidean space. Then
$\operatorname{ConvexHull} : \mathcal{A}(\mathbf{X}) \to
\mathcal{A}(\mathbf{X})$ is computable.
\end{prop}
\proof
Let $(K_n)_{n \in \mathbb{N}}$ be an effective enumeration of the convex hulls generated by finitely many points with rational coordinates. The interiors are uniformly computable, too. Moreover, in a computably compact space, closed sets are uniformly compact, inclusion of compact sets in open sets is recognizable, and countable intersection is computable on closed sets. Thus, the following equation determines computability of $\operatorname{ConvexHull}$:
\[\operatorname{ConvexHull}(A) = \bigcap_{n \in \{i \mid A \subseteq K_i^\circ\}} K_n\eqno{\qEd}\]

\begin{defi}[Closed Choice \cite{brattka3}]
Let $\mathbf{X}$ be a represented space. Then the {\em closed choice} operation
$\C_{\mathbf{X}}:\subseteq \mathcal{A}(\mathbf{X})\mto \mathbf{X}$ of this space is defined by
$x \in \C_\mathbf{X}(A)$ iff $x \in A$,
with $\dom(\C_X):=\{A\in \mathcal{A}(\mathbf{X}):A\not=\emptyset\}$.
\end{defi}

Intuitively, $\C_\mathbf{X}$ takes as input a non-empty closed set in negative representation (i.e.\ given by the capability to recognize the complement)
and produces an arbitrary point of this set as output.

\begin{defi}
For a represented space $\mathbf{X}$ and $1 \leq n \in \mathbb{N}$, let $\C_{\mathbf{X},\sharp=n} := \C_{\mathbf{X}}|_{\{A \in \mathcal{A}(\mathbf{X}) \mid |A| = n\}}$ and $\C_{\mathbf{X},\sharp\leq n} := \C_{\mathbf{X}}|_{\{A \in \mathcal{A}(\mathbf{X}) \mid 1 \leq |A| \leq n\}}$.
\end{defi}
More generally, for any choice principle the subscript $\sharp = n$ denotes the restriction to sets of cardinality $n$, and the subscript $\sharp \leq n$ to non-empty sets of cardinality less or equal than $n$. In the same spirit, the subscript $\lambda > \epsilon$ denotes the restriction to sets of outer radius greater than $\epsilon$, and $\mu > \epsilon$ the restriction to those sets where some value $\mu$ is greater than $\varepsilon$.

\begin{defi}
Let $\xc{n} := \C_{[0, 1]^n }|_{ \{A \in \mathcal{A}([0, 1]^n) \mid A \textnormal{ is convex }\}}$.
\end{defi}

The proof of the following proposition has been inspired by the proof of \cite[Theorem 3.1]{kreinovich} by \name{Longpr\'e} \etal, which the proposition generalizes in some sort. In fact, the study of $\C_{\uint, \sharp=m}$ is quite closely related to the theme of \cite{kreinovich}.
\begin{prop}
\label{iccccms}
Let $\mathbf{X}$ be a computably compact computable metric space. Then $\C_{\mathbf{X}, \sharp=n} \leq_{sW} \C_{\Cantor, \sharp=n}$ and $\C_{\mathbf{X}, \sharp \leq n} \leq_{sW} \C_{\Cantor, \sharp \leq n}$.
\begin{proof}
We associate a labeled binary infinite tree $T_\mathbf{X}$ with the space $\mathbf{X}$, where the vertices of some layers are labeled by closed balls in $\mathbf{X}$.

The root is labeled by $\mathbf{X}$. Then we find a finite open cover of $\mathbf{X}$ by open balls $B(x_1, 2^{-1}), \ldots, B(x_{2^k}, 2^{-1})$ using the computable dense sequence and the computable compactness provided by $\mathbf{X}$. We label the $k$-th layer of the tree with the closed balls $\overline{B}(x_i, 2^{-1})$. For the next step, each closed ball $\overline{B}(x_i, 2^{-1})$ (which is computably compact as a computably closed subset of a computably compact space) is covered by finitely many $B(x_{i,j}, 2^{-2})$ (which we can find by computable compactness of $\overline{B}(x_i, 2^{-1})$), and we then use $\overline{B}(x_{i,j}, 2^{-2})$ as labels for a suitable layer further down the tree.

In the next step, we cover each $\overline{B}(x_i,2^{-1}) \cap \overline{B}(x_{i,j},2^{-2})$ by finitely many open balls of radius $2^{-3}$, and so on. This process is iterated indefinitely, yielding finer and finer coverings of the space at each layer. For any labeled vertex we keep testing whether the intersection of the closed balls on the path leading to it is empty (which is computable by computable compactness), and if this is detected, the corresponding subtree is pruned at its current depth. Now any infinite path through the tree computably determines a point in $\mathbf{X}$ obtained by taking the intersection of all closed balls occurring as labels on the paths; and any point in $\mathbf{X}$ can be obtained from some infinite path.

Any closed subset of a computably compact space is compact (in a uniform way), so we can assume the input to $\C_{\mathbf{X}, \sharp=n}$ (resp.~$\C_{\mathbf{X}, \sharp\leq n}$) to be a compact set $A$ of cardinality $n$ (resp. less-or-equal to $n$). On any labeled layer of the tree beyond the $\left (1 + \log n\right )$-th, there are $n$ vertices such that the union of the intersections of the labels on the paths leading to them covers $A$. It is recognizable when an open set includes a compact set, so we will find suitable $n$ vertices eventually. Also, we can require that the vertices chosen on one level are actually below those chosen on the previous level. On the unlabeled layers of the tree, we simply remove all vertices that do not have a remaining labeled vertex beneath them. If we ever recognize that a proper subset of the $n$ vertices at some level already covered $A$, we prune the corresponding subtrees of the unneeded vertices at their current depth. With this process, we compute a subtree $T_A$ of $T_X$.

As $T_A$ has no more than $n$ vertices per layer, it is a name for a closed subset of $\Cantor$ with no more than $n$ points. If $A$ has exactly $n$ points, then from some layer onwards, $T_A$ will be the union of $n$ distinct pathes, hence be a name for a closed subset of $\Cantor$ with exactly $n$ points. As mentioned above, any infinite path through $T_A$ induces a point in $\mathbf{X}$, which as $A$ is closed, will actually fall in $A$.
\end{proof}
\end{prop}

It is rather obvious that if $\mathbf{X}$ is a co-c.e. closed subspace of $\mathbf{Y}$, then $\C_{\mathbf{X}, \sharp = n} \leq_{sW} \C_{\mathbf{Y}, \sharp = n}$ and $\C_{\mathbf{X}, \sharp \leq n} \leq_{sW} \C_{\mathbf{Y}, \sharp \leq n}$ (compare \cite[Section 4]{paulybrattka}). We recall that a computable metric space $\mathbf{X}$ is called \emph{rich}, if it has a subspace that is computably isomorphic to Cantor space (then this subspace automatically is co-c.e.~closed). \cite[Proposition 6.2]{gherardi3} states that any non-empty computable metric space without isolated points is rich.

\begin{cor}
\label{corr:icnrichcompact}
Let $\mathbf{X}$ be a rich computably compact computable metric space. Then $\C_{\mathbf{X}, \sharp = n} \equiv_{sW} \C_{\Cantor, \sharp = n}$ and $\C_{\mathbf{X}, \sharp \leq n} \equiv_{sW} \C_{\Cantor, \sharp \leq n}$.
\end{cor}

By inspection of the proof of Proposition \ref{iccccms}, we notice that the names produced there as inputs to $\C_{\Cantor, \sharp=n}$ or $\C_{\Cantor, \sharp \leq n}$ have a specific form: If we consider the closed subsets of Cantor space to be represented as the sets of infinite paths of infinite binary trees, the trees involved will have exactly $n$ vertices on all layers admitting at least $n$ vertices in a complete binary tree. The names used for $\C_{\Cantor, \sharp=n}$ moreover have the property that from some finite depths onwards, all vertices have exactly one child. We shall denote by $\ic{\leq n}$ (by $\ic{=n}$) the problem of finding an infinite path through a tree having exactly $n$ vertices from the $\lceil \log n \rceil$-th layer onwards (and where eventually each vertex has exactly one child)\footnote{Note that $\ic{=n}$ and $\ic{\leq n}$ are not restrictions of $\C_{\Cantor, \sharp=n}$ or $\C_{\Cantor, \sharp \leq n}$, but the realizers of the former problems are restrictions of the realizers of the latter.}. We directly conclude $\ic{=n} \equiv_{sW} \C_{\Cantor, \sharp=n} \equiv_{sW} \C_{\uint^k, \sharp=n}$ and $\ic{\leq n} \equiv_{sW} \C_{\Cantor, \sharp \leq n} \equiv_{sW} \C_{\uint^k, \sharp \leq n}$ for $k > 0$.

\section{Relative separation techniques}
\label{sec:relativetechniques}
The relative separation techniques to be developed in this section do not enable us to prove separation results just on their own; instead they constitute statements that some reduction $f \leqW g$ implies some reduction $f' \leqW g'$ , so by contraposition $f' \nleqW g'$ (which may be easier to prove) implies $f \nleqW g$. A particular form of these implications are absorption theorems. These show that for special degrees $h$, whenever $f$ has a certain property, then $f \leqW g \star h$ (or $f \leqW h \star g$) implies $f \leqW g$. A known result of this form is the following:

\begin{thm}[{\name{Brattka}, \name{de Brecht} \& \name{Pauly} \cite[Theorem 5.1]{paulybrattka}\footnote{The precise statement of \cite[Theorem 5.1]{paulybrattka} is weaker than the one given here, but a small modification of the proof suffices to obtain the present form. The only property of computable metric spaces used in that proof is that from a compact singleton $\{y\}$ the point $y$ could have been extracted. This, however, is just the definition of computable admissibility. Moreover, replacing the parallel product with the sequential one has no significant impact on the structure of the proof.}}]
Let $\mathbf{X}$, $\mathbf{Y}$ be represented spaces, and $\mathbf{Y}$ be computably admissible (cf.~\cite{pauly-synthetic-arxiv}, following \name{Schr\"oder} \cite{schroder}). Let $f : \mathbf{X} \to \mathbf{Y}$ be single-valued. Then $f \leqW \C_{\Cantor} \star g$ implies $f \leqW g$.
\end{thm}

We call a Weihrauch-degree a \emph{fractal}, if each of its parts is again the whole. The concept was introduced by \name{Brattka}, \name{de Brecht} and \name{Pauly} in \cite{paulybrattka} as a criterion for a degree to be join-irreducible (all fractals are join-irreducible, cf.~Lemma \ref{lemma:fractalsjoinirreducible}).

\begin{defi}
\label{def:fractal}
We call $f : \mathbf{X} \mto \mathbf{Y}$ a \emph{fractal}, iff there is some $g :\subseteq \Baire \mto \mathbf{Z}$, $f \equivW g$ such that for any clopen $A \subseteq \Baire$, either $g|_A \equiv_{W} f$ or\footnote{Note that $g|_A \equivW 0$ happens if and only if $A \cap \dom(g) = \emptyset$.} $g|_A \equivW 0$. If we can choose $g$ to be total, we call $f$ a \emph{closed fractal}.
\end{defi}

We will prove two absorption theorems, one for fractals and one for closed fractals. These essentially state that certain Weihrauch degrees are useless in solving a (closed) fractal.

\begin{thm}[Fractal absorption]
\label{theo:fractalabsorption}
If $f$ is a fractal, then $f \leqW g \star \C_{\{1, \ldots, n\}}$ implies $f \leqW g$.
\begin{proof}
We prove that $f \leqW g \star \C_{\{1, \ldots, n\}}$ implies $f \leqW g \star \C_{\{1, \ldots, n - 1\}}$ for fractal $f$ and $n > 1$, and then iteration together with $h \star \C_{\{1\}} \equivW h$ does the rest. We make a case distinction for this: First, assume that the reduction $f \leqW g \star \C_{\{1, \ldots, n\}}$ always uses the input $\{1, \ldots, n\}$ for $\C_{\{1, \ldots, n\}}$. Then replacing $\C_{\{1, \ldots, n\}}$ by the constant computable function $1$ works equally well, and we get $f \leqW g$ directly. Otherwise, there is some input $p$ for $f$, such that some $i \in \{1, \ldots, n\}$ is not contained in the input used for $\C_{\{1, \ldots, n\}}$. But then, $i$ has to be removed at some finite stage, when only a finite prefix $p_{\leq k}$ has been read. Restricting $f$ to those inputs starting with $p_{\leq k}$ does not change its Weihrauch degree (as $f$ is a fractal, and $f_{p_{\leq k}\Cantor} \equivW 0$ cannot happen, as $\dom(f_{p_{\leq k}\Cantor}) \neq \emptyset$). But then, $i$ is never contained in the set used as input for $\C_{\{1, \ldots, n\}}$, hence, $\C_{\{1, \ldots, n - 1\}}$ can be used instead (after exchanging $i$ and $n$).
\end{proof}
\end{thm}

\subsection{Baire Category Theorem as separation technique}
The absorption theorem for closed fractals is a consequence of the Baire Category Theorem, and was first employed as a special case in \cite[Proposition 4.9]{brattka3} by \name{Brattka} and \name{Gherardi}.

\begin{thm}[Closed fractal absorption]
\label{theo:closedfractalabsorption}
If $f$ is a closed fractal, then $f \leqW g \star \C_\mathbb{N}$ implies $f \leqW g$.
\begin{proof}
The degree $g \star \C_\mathbb{N}$ has a representative of the form $g' \circ \left (\id_\Baire \times \C_\mathbb{N} \right)$ with $g' \equivW g$, as shown in \cite{paulybrattka4}. W.l.o.g. assume that $f$ witnesses its own closed fractality. Let the inner reduction witness for $f \leqW  g' \circ \left (\id_\Baire \times \C_\mathbb{N} \right)$ be of the form $\langle H_1, H_2\rangle$. In particular, $H_2 : \Baire \to \dom(\C_\mathbb{N}\psi_-^\mathbb{N})$ is a computable map.

The closed sets $A_n = \{p \mid n \in \psi_-^\mathbb{N}(p)\}$ cover $\dom(\C_\mathbb{N}\circ \psi_-^\mathbb{N}) \subseteq \Baire$, and the corresponding restrictions $(\C_\mathbb{N})_{A_n}$ are computable for each $n \in \mathbb{N}$ by virtue of the constant function with value $n$ being a suitable realizer.  The closed sets $H_2^{-1}(A_n)$ cover $\dom(f) = \Baire$. Thus, we can apply the Baire Category Theorem, and conclude that there exists some $n_0$ such that $H_2^{-1}(A_{n_0})$ contains some non-empty clopen ball. As $f$ is a fractal, we know: $$f \leqW f_{H_2^{-1}(A_{n_0})} \leqW (g' \star \left (\id_\Baire \times (\C_\mathbb{N})_{A_{n_0}}\right) \leqW g' \equivW g$$
\end{proof}
\end{thm}

The preceding result occasionally is more useful in a variant adapted directly to choice principles in the r\^ole of $g$. For this, we recall the represented space $\mathbb{R}_{>}$, in which decreasing sequences of rational numbers are used to represent their limits as real numbers. We use $\overline{\mathbb{R}}_{>}$ to denote $\mathbb{R}_> \cup \{+\infty\}$, where $+\infty$ is represented by an empty sequence of rationals. Note that $\id : \mathbb{R} \to \mathbb{R}_{>}$ is computable but lacks a computable inverse. A generalized measure\footnote{As demonstrated in \cite{schroder2} (see also \cite{collins4}), one can obtain a canonical representation of the space of probability measures on some space $\mathbf{X}$ by restricting $\mathcal{C}(\mathcal{O}(\mathbf{X}), \mathbb{R}_{<})$ to those functions satisfying the properties of probability measures. By moving to the complement, one arrives at the present setting.} on some space $\mathbf{X}$ is a continuous function $\mu : \mathcal{A}(\mathbf{X}) \to \mathbb{R}_>$ taking only non-negative values. The two variants of the Baire Category theorem as separation technique are connected by the following result:

\begin{prop}
Define $\textsc{Lb} : \{x \in \mathbb{R}_{>} \mid x > 0\} \to \mathbb{N}$ via $\textsc{Lb}(x) = \min \{n \in \mathbb{N} \setminus \{0\} \mid n^{-1} \leq x\}$. Then $\textsc{Lb} \equiv_{sW} \C_\mathbb{N}$.
\begin{proof}
Given a $\rho_{<}$-name of $x$, the property $n^{-1} \leq x$ is refutable: If $x < n^{-1}$, then the rational sequence approaching $x$ from above must pass $n^{-1}$ at some point. Hence we can compute $\{n \in \mathbb{N} \setminus \{0\} \mid x < n^{-1}\} \in \mathcal{O}(\mathbb{N})$. Finding the maximum in this set is (strongly) reducible to $\C_\mathbb{N}$ (e.g. by \cite[Theorem 4.3.1.24]{paulyphd}), it remains to increment it by $1$.

For the other direction we present a reduction from $\C_\mathbb{N}$. Once all integers from $0$ to $k$ have been encountered in the input to $\C_\mathbb{N}$, we print the rational $(k + 1.5)^{-1}$ (with sufficiently many repetitions to ensure an infinite output). If $n_0$ is the smallest solution to $\C_\mathbb{N}$, this produces a $\rho_{>}$-name of $(n_{0} + 0.5)^{-1}$, hence application of $\textsc{Lb}$ will return $n_0 + 1$.
\end{proof}
\end{prop}

The preceding result indirectly shows how a closed choice principle for some class $\mathfrak{A} \subseteq \mathcal{A}(\mathbf{X})$ of closed sets with positive generalized measure $\mu$ can be decomposed into the slices with fixed lower bounds $\mu > n^{-1}$. For this, we recall the infinitary coproduct (i.e.~disjoint union) $\coprod_{n \in \mathbb{N}}$ defined both for represented spaces and multivalued functions between them via $\left ( \coprod_{n \in \mathbb{N}} f_n \right )(i, x) = (i, f_i(x))$.

\begin{cor}
\label{corr:boundingmeasures}
$\mathrm{C}_{\mathbf{X}}|_{\mathfrak{A},\mu > 0} \leqW \left ( \coprod_{n \in \mathbb{N}} \mathrm{C}_{\mathbf{X}}|_{\mathfrak{A},\mu > n^{-1}} \right ) \star \C_\mathbb{N}$
\end{cor}

\begin{lem}[{$\sigma$-join irreducibility of fractals \cite[Lemma 5.5]{paulybrattka}}]
\label{lemma:fractalsjoinirreducible}
Let $f$ be a fractal and satisfy $f \leqW \coprod_{n \in \mathbb{N}} g_n$. Then there is some $n_0 \in \mathbb{N}$ such that $f \leqW g_{n_0}$.
\end{lem}

\begin{thm}
\label{theo:closedfractalabsorption2}
Let $f$ be a closed fractal such that $f \leqW \C_{\mathbf{X}}|_{\mathfrak{A},\mu>0}$. Then there is some $n \in \mathbb{N}$ such that $f \leqW \C_{\mathbf{X}}|_{\mathfrak{A},\mu>n^{-1}}$.
\begin{proof}
From Corollary \ref{corr:boundingmeasures} we deduce that $f \leqW \left ( \coprod_{n \in \mathbb{N}} \mathrm{C}_{\mathbf{X}}|_{\mathfrak{A},\mu > n^{-1}} \right ) \star \C_\mathbb{N}$. Then Theorem \ref{theo:closedfractalabsorption} implies $f \leqW \left ( \coprod_{n \in \mathbb{N}} \mathrm{C}_{\mathbf{X}}|_{\mathfrak{A},\mu > n^{-1}} \right )$. By Lemma \ref{lemma:fractalsjoinirreducible} there has to be some $n_0$ with $f \leqW \mathrm{C}_{\mathbf{X}}|_{\mathfrak{A},\mu > n_0^{-1}}$.
\end{proof}
\end{thm}

Before ending this subsection, we shall provide some useful examples of generalized measures that are not already measures. Note that we only use the implication $2. \rightarrow 1.$ from Proposition \ref{prop:lambda}, while Proposition \ref{prop:iota} is only included for completeness. For some subset $A \subseteq \mathbf{X}$ of a metric space, the outer radius $\lambda$ is defined via $\lambda(A) = \inf \{r \geq 0 \mid \exists x \ A \subseteq B(x, r)\}$, and the inner radius $\iota$ is defined via $\iota(A) = \sup \{r \geq 0 \mid \exists x \ B(x, r) \subseteq A\}$.
\begin{prop}
\label{prop:lambda}
For a computable metric space $\mathbf{X}$, the following are equivalent:
\begin{enumerate}
\item $\lambda : \mathcal{A}(\mathbf{X}) \to \overline{\mathbb{R}}_>$ is computable.
\item $\mathbf{X}$ is computably compact.
\end{enumerate}
\begin{proof}\hfill
\begin{description}
\item[$1. \Rightarrow 2.$] If $\mathbf{X}$ is a computable singleton, then it is computably compact anyway. If not, there are two distinct computable points $x, y \in \mathbf{X}$. From some $A \in \mathcal{A}(\mathbf{X})$, we can compute $A \cup \{x\} \in \mathcal{A}(\mathbf{X})$ and $A \cup \{y\} \in \mathcal{A}(\mathbf{X})$. Furthermore, $\max : \overline{\mathbb{R}}_> \times \overline{\mathbb{R}}_> \to \overline{\mathbb{R}}_>$ and $\mathalpha{>} : \mathbb{R} \times \overline{\mathbb{R}}_> \to \mathbb{S}$ are computable, so assuming that $\lambda$ is computable, we can compute $A \mapsto \left [0.1d(x, y) > \max (\lambda(A \cup \{x\}), \lambda(A \cup \{y\}))\right ]$. Now this expression will evaluate to \emph{true} if and only if $A$ is empty: If $A =\emptyset$, the righthand side is zero but the lefthand side is not; if $z \in A$ then $\lambda(A\cup \{x\}) \geq \frac{1}{2}d(x,z)$, and the triangle inequality can be invoked to arrive at the contradiction $0.4d(x,y) > d(x,y)$. Thus, we have demonstrated that $\textsc{IsEmpty} : \mathcal{A}(\mathbf{X}) \to \mathbb{S}$ is computable. This in turn is the definition of computable compactness.
\item[$2. \Rightarrow 1.$] In a computably compact space, closed sets are compact (in a uniform way). This in turn makes $A \subseteq B(x, r)$ recognizable. Next, we point out that in $\lambda(A) = \inf \{r \mid \exists x \ A \subseteq B(x, r)\}$ it suffices to have $x$ range over the dense basic sequence in $\mathbf{X}$, and $r$ over the positive rationals. Finally, $\inf : \mathcal{O}(\mathbb{Q}^+) \to \overline{\mathbb{R}}_>$ is computable, hence the claim follows from the definition.\qedhere
\end{description}
\end{proof}
\end{prop}

\begin{prop}
\label{prop:iota}
Let $\mathbf{X}$ be a computably compact computable metric space. Then $\iota : \mathcal{A}(\mathbf{X}) \to  \overline{\mathbb{R}}_>$ is computable.
\begin{proof}
Let $(q_i)_{i \in \mathbb{N}}$ be a computable dense sequence in $\mathbf{X}$. The compactness of $\mathbf{X}$ allows us to enumerate all tuples $\langle r, i_1, \ldots, i_n\rangle$ such that $\mathbf{X} \subseteq \bigcup_{j \leq n} B(q_{i_j}, r)$. Given some closed set $A \in \mathcal{A}(\mathbf{X})$, we can narrow this down to those tuples where additionally $\forall j \leq n \ q_{i_j} \notin A$. Now note that $\iota(A)$ is the infimum of all $r$ occurring in this enumeration.

To substantiate the latter claim, we observe the following: If $r > \iota(A)$, then for any $x \in A$ the set $B(x,r) \cap A^C$ is non-empty. As this is an open set, it will then contain some basic point $q_x$. From $A \subseteq \bigcup_{x \in A} B(q_x,r)$ and compactness of $A$ we see that some finite number of basis points is sufficient to witness non-emptiness of $B(x,r) \cap A^C$ for all $x \in A$ simultaneously. We can safely add finitely more points to also cover the rest of $\mathbf{X}$. The procedure above will eventually find such a collection, hence $r$ is taken into account for the infimum, and we see that we cannot compute too large a value.

For the converse direction, let us assume that $\mathbf{X} \subseteq \bigcup_{j \leq n} B(x_j, r)$ for some points $x_j \notin A$. If there were some ball $B(x,r) \subseteq A$, we would arrive at a contradiction as follows: As $x \in \mathbf{X}$, there is some $x_j \notin A$ with $x \in B(x_j, r)$. But by symmetry, then also $x_j \in B(x,r) \subseteq A$. Thus, we conclude that any $r$ occurring in our enumeration actually is an upper bound for $\iota(A)$, hence the computation works correctly.
\end{proof}
\end{prop}
\subsection{Large radius technique}
Given a closed fractal $f$, Theorem \ref{theo:closedfractalabsorption2} allows us to bound away from 0 any
positive generalized measure on the closed sets that are used to compute
the function $f$. The separation technique to be developed next bounds
away from 0 only a specific generalized measure -- the outer radius --
yet requires neither positivity nor the closed fractal property.

For a computable metric space $\mathbf{X}$, $\varepsilon > 0$ and some class $\mathfrak{A} \subseteq \mathcal{A}(\mathbf{X})$, we introduce: \[X_\varepsilon(\mathfrak{A}) = \overline{\psi_-^{-1}\left(\{A \in \mathfrak{A} \mid \forall x \in \mathbf{X} \exists B \in \mathfrak{A} \ B \subseteq A \setminus B(x, \varepsilon)\}\right)} \subseteq \Baire\] This means that the names in $X_\varepsilon(\mathfrak{A})$ are for sets large enough such that arbitrarily late an arbitrary ball of radius $\varepsilon$ can be removed from them, and still a closed set in the class $\mathfrak{A}$ remains as a subset. As $X_\varepsilon(\mathfrak{A})$ is a set of names for instances of $\C_\mathbf{X}$, rather than a set of instances itself, we have to use the generalized restriction $f_A$ introduced on page \pageref{page:fa} and study $(\C_{\mathbf{X}})_{X_\varepsilon(\mathfrak{A})} = (\C_{\mathbf{X}}|_{\mathfrak{A}})_{X_\varepsilon(\mathfrak{A})}$ rather than the meaningless $\C_{\mathbf{X}}|_{X_\varepsilon(\mathfrak{A})}$.

\begin{question}
\label{question:fractals}
If $\C_{\mathbf{X} }|_{\mathfrak{A}}$ is a fractal, is $(\C_{\mathbf{X}})_{X_\varepsilon(\mathfrak{A})}$ too\footnote{Clearly, if this is true, then whenever $\C_{\mathbf{X} }|_{\mathfrak{A}}$ is a closed fractal, so is $(\C_{\mathbf{X}})_{X_\varepsilon(\mathfrak{A})}$.}?
\end{question}

We proceed to show that a reduction between choice principles has to map sets large in this sense to sets with large outer radius (denoted by $\lambda$).

\begin{lem}[Large Radius Principle]
\label{lemma:largediameter}
Let $H$ and $K$ witness a reduction \linebreak $\C_{\mathbf{X}}|_{\mathfrak{A}} \leqW \C_{\mathbf{Y}}|_{\mathfrak{B}}$, where $\mathbf{Y}$ is compact and $\mathfrak{A} \subseteq \mathcal{A}(\mathbf{X})$, $\mathfrak{B} \subseteq \mathcal{A}(\mathbf{Y})$. Then \[\forall p \in \dom(\C_{\mathbf{X}}|_{\mathfrak{A}}\psi_-^\mathbf{X}) \ \forall \varepsilon > 0 \ \exists n \in \mathbb{N} \ \exists \delta > 0 \ \forall q \ \left ( q \in X_\varepsilon(\mathfrak{A}) \cap B(p, 2^{-n}) \Rightarrow \lambda\psi^\mathbf{Y}_- H (q) > \delta\right )\]
\end{lem}
\begin{proof}
Assume the claim were false, and let $p \in \dom(\C_{\mathbf{X}}|_{\mathfrak{A}}\psi^\mathbf{X}_-)$ and $\varepsilon > 0$ be witness for the negation. There has to be a sequence $(p_n)_{n \in \mathbb{N}}$ such that $p_n \in X_\varepsilon(\mathfrak{A})$, $d(p, p_n) < 2^{-n}$ and $\lambda\psi^\mathbf{Y}_- H(p_n) < 2^{-n}$. As the $p_n$ converge to $p$ and $H$ is continuous, we conclude that $\lim_{n \to \infty} H(p_n) = H(p)$. For the closed sets represented by these sequences, this implies $\left ( \bigcap_{n \in \mathbb{N}} \overline{\bigcup_{i \geq n} \psi_-^\mathbf{Y} H (p_i)} \right ) \subseteq \psi_-^\mathbf{Y} H(p)$. As $\mathbf{Y}$ is compact, the left hand side contains some point $x$.

As $x \in \psi_-^\mathbf{Y} H (p)$, for any $q \in \delta_\mathbf{Y}^{-1}(\{x\})$ it is the case that $\langle p, q\rangle \in \dom(K)$. We fix such a $q$ and $y = \delta_\mathbf{X}(K(\langle p, q\rangle))$. By continuity, there is some $N \in \mathbb{N}$ such that for any $\langle p', q'\rangle \in (B(p, 2^{-N})\times  B(q, 2^{-N})) \cap \dom(\delta_\mathbf{X}K)$ it follows that $\delta_\mathbf{X}K(\langle p', q'\rangle) \in B(y, \varepsilon)$.

By choice of $x$, for any $i \in \mathbb{N}$ there is some $k_i \geq i$ such that $d(x, \psi_-^\mathbf{Y} H (p_{k_i})) < 2^{-i}$. By choice of the $p_n$, this in turn implies $\psi_-^\mathbf{Y} H (p_{k_i}) \subseteq B(x, 2^{-i} + 2^{-k_i+1})$. Let $I \in \mathbb{N}$ be large enough, such that for any $x' \in B(x, 2^{-I} + 2^{-k_I+1})$ it follows that $\delta_\mathbf{Y}^{-1}(x') \cap B(q, 2^{-N}) \neq \emptyset$. The inclusion $\psi_-^\mathbf{Y} H (p_{k_I}) \subseteq B(x, 2^{-I} + 2^{-k_I+1})$ of a compact set in an open set implies that there is some $L > k_I$ such that for all $p' \in B(p_{k_I}, 2^{-L}) \cap \dom(\C_{\mathbf{X}}|_{\mathfrak{A}}\psi_-^\mathbf{X})$ it holds that $\psi_-^\mathbf{Y} H p' \subseteq B(x, 2^{-I} + 2^{-k_I+1})$.

The choice of $p_{k_I}$, $L$ and the point $y \in \mathbf{X}$ ensures that our reduction may answer any valid input to $\C_{\mathbf{X}}|_{\mathfrak{A}}$ sharing a prefix of length $L$ with $p_{k_I}$ with a name of some point $y' \in B(y, \varepsilon)$. However, as we have $p_{k_I} \in X_\varepsilon(\mathfrak{A})$, we can extend any long prefix of $p_{k_I}$ to a name of a set not intersecting the ball $B(y, \varepsilon)$ -- this means, our reduction would answer incorrectly, and we have found the desired contradiction.
\end{proof}

\begin{cor}[Large Radius Principle for fractals]
\label{corr:ldpfractals}
Let $\C_{\mathbf{X}}|_{\mathfrak{A}}$ be a fractal, $\mathbf{Y}$ be compact and $\C_{\mathbf{X}}|_{\mathfrak{A}} \leqW \C_{\mathbf{Y}}|_{\mathfrak{B}}$. Then for every $\varepsilon > 0$ there exists a $\delta > 0$ such that \[(\C_{\mathbf{X}})_{X_\varepsilon(\mathfrak{A})} \leqW \C_{\mathbf{Y}}|_{\mathfrak{B}, \lambda > \delta }\]
\begin{proof}
Let $g : \subseteq \Baire \mto \mathbf{Z}$ witness fractality of $\C_{\mathbf{X}}|_{\mathfrak{A}}$. Let $K, H$ witness the reduction $\C_{\mathbf{X}}|_{\mathfrak{A}} \leqW \C_{\mathbf{Y}}|_{\mathfrak{B}}$ and let $K', H'$ witness the reduction $g \leqW \C_{\mathbf{X}}|_{\mathfrak{A}}$. Pick some $p' \in \dom(g)$. Choose $\delta > 0$ and $n \in \mathbb{N}$ satisfying the property in Lemma \ref{lemma:largediameter} for $\varepsilon$ and $p := H'(p')$. Continuity of $H'$ implies the existence of some $n' \in \mathbb{N}$ such that $H'[B(p', 2^{-n'})] \subseteq B(p, 2^{-n})$. By Lemma \ref{lemma:largediameter}, $\langle x, y\rangle \mapsto K'(\langle x, K(\langle H'(x), y\rangle)\rangle)$ and $H \circ H'$ witness $g|_{B(p', 2^{-n'})} \leqW \C_{\mathbf{Y}}|_{\mathfrak{B}, \lambda > \delta }$. Fractality provides $(\C_{\mathbf{X}})_{X_\varepsilon(\mathfrak{A})} \leqW g|_{B(p', 2^{-n'})}$, thus the claim is demonstrated.
\end{proof}
\end{cor}

\section{Separation results for finite and convex choice}

\begin{figure}[htbp]
\begin{tikzpicture}
  \matrix (m) [matrix of math nodes, row sep=3em, column sep=3em]
    { \phantom{\xc{n}} & \boxed{\xc{1}}  & \xc{2}  & \xc{n}  & \phantom{\xc{n}} & \C_{\uint} \\
    & \ic{\leq2}  & \ic{\leq3}  & \ic{\leq n+1}  &  \phantom{\ic{\leq n+1}} & \\
      & \ic{=2}  & \ic{=3}  & \ic{=n+1}  & \phantom{\ic{=n+1}} & \C_\mathbb{N}\\
      1 \equivW \C_{\{0\}} & \C_{\{0, 1\}} & \C_{\{0, 1, 2\}} & \C_{\{0, \ldots, n\}} & \phantom{\C_{\{0, 1\}}} &  \boxed{\xc{1}}\\ };
  { [start chain] \chainin (m-1-1);
      \chainin (m-1-2);
    \chainin (m-1-3);
    \chainin (m-1-4);
    \chainin (m-1-5);
        \chainin (m-1-6);}
    { [start chain] \chainin (m-2-2);
    \chainin (m-2-3);
    \chainin (m-2-4);
    \chainin (m-2-5); }
    { [start chain] \chainin (m-3-2);
    \chainin (m-3-3);
    \chainin (m-3-4);
    \chainin (m-3-5);
        \chainin (m-3-6);  }
      { [start chain] \chainin (m-4-1);
    \chainin (m-4-2);
      { [start branch=A] \chainin (m-3-2); \chainin (m-2-2); \chainin (m-1-2);}
    \chainin (m-4-3);
      { [start branch=B] \chainin (m-3-3); \chainin (m-2-3); \chainin (m-1-3);}
    \chainin (m-4-4);
           { [start branch=C] \chainin (m-3-4); \chainin (m-2-4); \chainin (m-1-4);}
    \chainin (m-4-5);
          { [start branch=D] \chainin (m-3-5); \chainin (m-2-5); \chainin (m-1-5);}
    \chainin (m-4-6); }
\end{tikzpicture}
\caption{The reducibilities}
\label{figure:reducibilities}
\end{figure}
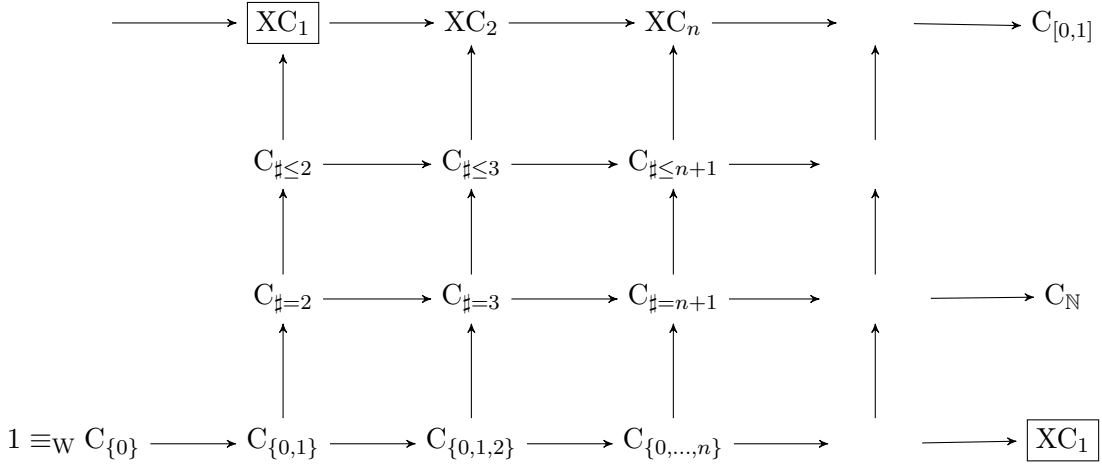

We now have the tools available to completely characterize the valid reductions between $\C_{\{0, \ldots, n\}}$, $\xc{m}$, $\ic{\leq i}$ and $\ic{=j}$. Figure \ref{figure:reducibilities} provides an overview -- the absence of an arrow (up to transitivity) indicates a proof of irreducibility. Two important results have already been established in the literature, namely $\C_{\{0,\ldots,n\}} \leW \C_{\{0,\ldots,n+1\}}$ in \cite{weihrauchb} by Weihrauch, and $\left (\coprod_{n \in \mathbb{N}} \C_{\{0,\ldots,n\}} \right) \leW \xc{1}$ in \cite{paulybrattka3,paulybrattka3cie} by Brattka and the authors.

Besides an application of the general techniques of the preceding section, more specialized proof methods are employed, some with a rather combinatorial character, others based on the properties of simplices. We also exhibit a technique suitable to transfer results from the compact case to the locally compact case.

\begin{observation}
\label{obs:fractals}
$\ic{=n}$ is a fractal, and $\ic{\leq n}$ even is a closed fractal.
\begin{proof}
Let $\delta : \subseteq \Baire \to \mathcal{T}_2$ be some standard representation of the binary trees, in particular let $\delta$ satisfy that any tree identical on the first $n$ levels to some tree $\delta(p)$ has a $\delta$-name $q$ with $d(p, q) < 2^{-n}$. Now for any clopen $A \subseteq \Baire$ with $A \cap \dom(\ic{=n}\delta) \neq \emptyset$ (resp.~$A \cap \dom(\ic{\leq n}\delta) \neq \emptyset$) we see that $\ic{=n} \leqW \ic{=n}\delta|_A$ (resp.~$\ic{\leq n} \leqW \ic{\leq n}\delta|_A$), as we append the infinite input tree to some leaf of a finite tree chosen to ensure membership in $A$. The condition of having $n$ vertices on each level (besides the first $\lceil \log n\rceil$ ones) can easily be kept by pruning the remaining leafs of the finite tree in a delayed way. Any infinite path through the resulting tree provides an infinite path through the original tree as a tail.

To see that $\ic{\leq n}$ is even a closed fractal, we argue that $\dom(\ic{\leq n}\delta)$ is a computable retract of $\Baire$, i.e.~that there is a computable function $R : \Baire \to \dom(\ic{\leq n}\delta)$ with $R|_{\dom(\ic{\leq n}\delta)} = \id_{\dom(\ic{\leq n}\delta)}$. For this, notice that we can detect if some prefix $w$ cannot be extended to some $p \in \dom(\ic{\leq n}\delta)$, as this corresponds to the tree not having exactly $n$ vertices on some level beyond the $\left (\lceil \log n\rceil\right )$-th. Moreover, if some prefix $w$ can be extended to some $p \in \dom(\ic{\leq n}\delta)$, we can do so in a computable way (e.g.~by giving each leaf in the finite tree determined by $w$ exactly one successor in perpetuity). Now $\ic{\leq n} \circ \delta \circ R$ satisfies the criteria for $g$ in Definition \ref{def:fractal}.
\end{proof}
\end{observation}

\begin{observation}
\label{obs:fractals2}
$\xc{n}$ is a closed fractal.
\begin{proof}
We use a representation $\psi_-$ of the closed subsets $\mathcal{A}(\uint^n)$ with the property that for any finite word $w$ there is some non-degenerate rational hypercube $H$ such that any closed set $A \subseteq H$ has some $\psi_-$-name $p$ with $w \prec p$. Such a representation can be obtained from \cite[Proposition 3.4]{paulybrattka3}. Rescaling $\uint^n$ to $H$ and back now establishes $\xc{n} \leqW \xc{n}\circ \psi_-|_{w\Baire}$, i.e.~$\xc{n}$ is a fractal.

To see that $\xc{n}$ is a closed fractal, first note that the representation $\psi_-$ above can be chosen as total. By using the computable operator $\operatorname{ConvexHull}$ from Proposition \ref{prop:convexhull} we obtain a retract to the convex sets. Finally, as $\uint^n$ is compact, if we encounter a name of an empty set, we notice at some finite stage, and can modify the name accordingly. Hence, we even obtain a computable retract $R$ from all names of closed sets to names of non-empty convex closed sets. Now $\xc{n} \circ \psi_- \circ R$ witnesses that $\xc{n}$ is a closed fractal.
\end{proof}
\end{observation}

\begin{cor}
\label{corr:icnnleqc1n}
$\ic{=n} \nleqW \C_{\{0, \ldots, m\}}$ for all $n > 1$ and $m \in \mathbb{N}$.
\begin{proof}
Assume the reduction would hold for some $n, m \in \mathbb{N}$. Observation \ref{obs:fractals} allows us to use Theorem \ref{theo:fractalabsorption} to conclude $\ic{=n}$ to be computable - a contradiction for $n > 1$.
\end{proof}
\end{cor}

\begin{prop}[\footnote{We are grateful to a referee for suggestion the present simplified proof.}]
\label{prop:icnleqwcn}
$\ic{=n} \leqW \C_\mathbb{N}$
\begin{proof}
We show that $\ic{=n}$ is non-deterministically computable with advice space $\mathbb{N}$ and invoke \cite[Theorem 7.2]{paulybrattka}. We guess some $k \in \mathbb{N}$ such that beyond the $k$-th level in our tree each vertex has exactly one successor. A wrong guess can be detected and rejected eventually, while a correct guess allows us to compute an infinite path by simply extending in the unique possible way from some existing vertex on the $k$-th level onwards.
\end{proof}
\end{prop}

\begin{cor}
$\ic{\leq 2} \nleqW \ic{=n}$
\begin{proof}
Assume $\ic{\leq 2} \leqW \ic{=n}$ for some $n \in \mathbb{N}$. By Proposition \ref{prop:icnleqwcn}, this implies $\ic{\leq 2} \leqW \C_\mathbb{N}$. Observation \ref{obs:fractals} together with Theorem \ref{theo:closedfractalabsorption} would show $\ic{\leq 2}$ to be computable, contradiction.
\end{proof}
\end{cor}

\subsection{Combinatorial arguments}
\begin{prop}
\label{prop:cntoicnplusone}
$\C_{\{0, \ldots, n\}} <_{sW} \ic{=n+1}$.
\begin{proof}
We use $\C_{\uint,\sharp=n+1}$ in place of $\ic{=n+1}$; and employ the same property of the representation of the closed sets as in the proof of Observation \ref{obs:fractals2}. Fix $n + 1$ disjoint closed proper intervals in $\uint$. Start to produce a name for the closed set containing all the centers of the intervals. If any $i \in \{0, \ldots, n\}$ is removed from the input to $\C_{\{0, \ldots, n\}}$, all of the corresponding closed interval is removed from the input to $\ic{=n+1}$. The left-most remaining interval center has been approximated to some finite precision so far, hence there is still an open ball around it left in the current input to $\ic{=n+1}$. This ball is split into as many disjoint closed proper intervals as necessary to keep the cardinality condition.

Iterating this process produces a closed set containing exactly $n + 1$ points in the end, and any element is included in one of the initial intervals. As these are closed and disjoint, we can determine the index of an interval from a point. This constitutes a valid answer to the input for $\C_{\{0, \ldots, n\}}$.
\end{proof}
\end{prop}

\begin{prop}[Pigeonhole principle]
\label{prop:cntoicn}
$\C_{\{0, \ldots, n\}} \nleqW \ic{\leq n}$
\begin{proof}
Assume that $K$, $H$ would witness a reduction $\C_{\{0, \ldots, n\}} \leqW \C_{\Cantor,\sharp\leq n} \left (\equivW \ic{\leq n} \right )$. We consider their behaviour on an input $p$ representing the full set $\{0, \ldots, n\}$. $H$ will compute a name for some closed set $A \subseteq \Cantor$ consisting of the points $a_1, \ldots, a_k$ with $k \leq n$. We see that $\langle p, a_i\rangle \in \dom(K)$, and $K(\langle p, a_i\rangle) \in \{0, \ldots, n\}$. By the pigeonhole principle, there is some $J \in \{0, \ldots, n\}$ such that $J \neq K(\langle p, a_i\rangle)$ for all $i$.

By continuity of $K$, there is some $M_i$ such that $K(\langle q, b_i\rangle) = K(\langle p, a_i\rangle)$ for all $q, b_i$ with $d(p, q) < 2^{-M_i}$ and $d(a_i, b_i) < 2^{-M_i}$. Let $M := \max_{i \leq k} M_i$, and $A_M := \bigcup_{i \leq k} \{q \in \Cantor \mid d(q, a_i) \leq 2^{-M}\}$. Continuity of $H$ means that there is some $N \geq M$ such that for any $B \in (\psi_-\circ H)[p_{\leq N}\Cantor]$ it is the case that $B \subseteq A_M$. But with this, we have demonstrated that for any $q \in p_{\leq N}\Cantor$ as input to $\C_{\{0, \ldots, n\}}$ the reduction will eventually produce some $l \in \{0, \ldots, n\}$ with $l \neq J$. However, a name for $\{0, \ldots, n\}$ shares arbitrarily long prefixes with names for $\{J\}$, hence the reduction will fail.
\end{proof}
\end{prop}

\begin{cor}
$\ic{\leq n} \leW \ic{\leq n+1}$
\end{cor}

\begin{cor}
$\xc{1} \nleqW \ic{\leq n}$ for all $n \in\mathbb{N}$
\begin{proof}
By combining \cite[Proposition 7.2]{paulybrattka3} with \cite[Theorem
32]{paulyincomputabilitynashequilibria}, allows us to conclude that $\left (\coprod_{i \in \mathbb{N}} \C_{\{0,\ldots,i\}} \right) \leW \xc{1}$. In particular, $\xc{1} \leqW \ic{\leq n}$ would imply $\C_{\{0, \ldots, n\}} \leqW \ic{\leq n}$ and thus contradict Proposition \ref{prop:cntoicn}.
\end{proof}
\end{cor}

The following lemma serves to keep the algorithm employed in Proposition \ref{prop:activevertices} simple:
\begin{lem}
\label{lemma:activevertices-aux}
Restricting $\ic{=n}$ or $\ic{\leq n}$ to those trees where on each layer beyond the $\lceil \log n\rceil$-th at most one vertex has zero children (and hence at most one vertex has two children) does not change their Weihrauch-degree.
\begin{proof}
Only one direction of the reduction is non-trivial. The outer reduction witness $K$ is defined by $K(\langle p, q\rangle)(i) = q(in)$, i.e.~we take only every $n$-th bit of the infinite path through the derived tree in order to form the path through the original tree. The inner reduction witness expands any layer in the original tree to $n$ layers in the derived tree. If a vertex in the original tree has only the left (right) child, then in the derived tree, we add a tree of depth $n$ with only the left-most (right-most) branch in the derived tree. If a vertex in the original tree is the $k$-th vertex from the left to have two children, then we add the following tree of height $n$: The left-most branch until level $k$, then a split, followed by the right-most and the left-most branch only. If a vertex in the original tree is the $k$-th vertex from the left to have no children, we add a tree of height $k$ containing the left-most branch only (and ending in a leaf). Thus, the vertices with $0$ and $2$ children respectively from the original layer are spread out in a pairwise fashion, and the resulting tree satisfies the extra criterion.
\end{proof}
\end{lem}

\begin{prop}
\label{prop:activevertices}
$\ic{=n+1} \leq_{W} \ic{=2}^n$ and $\ic{\leq n+1} \leq_{W} \ic{\leq 2}^n$
\begin{proof}
We show how from a single infinite binary tree with $n + 1$ vertices per level beyond the $\lceil \log n+1\rceil$-th we can compute $n$ infinite binary trees with $2$ vertices per level beyond the first, such that knowing infinite paths through the latter trees allows us to pick an infinite path in the former. Moreover, the construction will ensure that if from some level onwards the original tree has exactly one successor per vertex, the same holds true for the derived trees. By Lemma \ref{lemma:activevertices-aux} we can assume freely that on each layer of the tree there is at most one leaf.

We shall call a vertex located at some layer less than $i$ \emph{active} at level $i$, if both its successors are extended by paths reaching the level $i$. If there are $n + 1$ vertices at each sufficiently large level, then there are $n$ active vertices a each sufficiently layer large $i$. The construction starts once the first few levels of the input tree have been seen, such that the number of vertices per level can reach $n + 1$ for the first time (it will remain at $n+1$ from there on). We place a distinct token from $\{1,\ldots,n\}$ on each of the active vertices.

The $n$ output trees (in the domain of $\ic{\leq 2}$ or even $\ic{=2}$) start off with the root and two children. The $k$-th output tree corresponds to the token $k$. As long as the vertex (in the original tree) holding the token $k$ remains active, we extend both paths of maximal length in the $k$-th output vertex.

Whenever the current layer in the original tree contains a leaf, then there is some active vertex with token $k$ that will cease to be active at the next level. This is the active vertex on the path to the leaf closest to the leaf. We take note whether it is the left or the right subtree of this active vertex that contains the leaf. By assumption, there is also a vertex on the same layer as the leaf that will have two children, i.e.~will become an active vertex. We move the token $k$ from the old to the new active vertex. In the $k$-th output tree, we cut the left or right subtree, depending on which subtree of the formerly active vertex contains the leaf, and give the last vertex of the other path two children.

It remains to describe how to find an infinite path through the input tree given infinite paths through the output trees. It is clear that the difficulty of finding a path extendable to an infinite one solely lies in the choice of which successor to pick at currently active vertices. Consider a vertex that became active at step $i$ and was labeled with the token $k$. Then by construction, choosing the same way as the path through the $k$-th output tree  at layer $i$ is safe.
\end{proof}
\end{prop}

\begin{exa}
Let us give an example with $n = 2$ as Figure \ref{fig:exampleactive} on Page \pageref{fig:exampleactive}. The upper half of the
table displays the step-by-step computation of two output trees (given an input
tree), and the lower half displays the step-by-step computation of an infinite
path in the input tree (given one infinite path in each output tree). In the
pictures we name only the branching vertices and the dead ends are marked by
solid circles. In the original tree, $a$, $c$ and $d$ are the vertices that are at some time active and carrying label $1$, whereas $b$ and $e$ are at some time active and carry label $2$.

The upper and lower parts of the table are shifted so that the
height of the input tree and the length of the computed path match in each
column. Also, in the lower part we give only the information about the paths in
the output trees that is relevant for the computation of the path in the input
tree.

\begin{figure}[htbp]
\begin{tabular}{|c|c|c|c|c|c|c|c|}
\hline
& Computing  & Computing &&&&& Receiving\\
&trees, & trees, &&&&& two paths.\\
&{\tiny generic} & {\tiny step 1.} & {\tiny step 2.} & {\tiny step 3.} & {\tiny step 4.} & {\tiny step 5.} &\\
&{\tiny start.} &&&&&&
\\\hline
{\tiny \begin{tabular}{c}\\Input tree\end{tabular}}&
\begin{tikzpicture}[baseline,level distance=5mm]
\node(root){{\tiny a}}[sibling distance=10mm]
        child{node{{\tiny b}}
                child{node{}}
                child{node{}}
        }
        child{node{} edge from parent[draw=none] child{node(right){} edge from
parent[draw=none]}};
                \draw   (root) -- (right)       ;
\end{tikzpicture}
&
\begin{tikzpicture}[baseline,level distance=5mm]
\node(root){\tiny a}[sibling distance=10mm]
        child{node{\tiny b}
                child{node{}
                        child{node{}}
                }
                child{node{}
                        child{node{}}
                }
        }
        child{node{} edge from parent[draw=none] child{node(right){} edge from
parent[draw=none] child{node{}}}};
                \draw   (root) -- (right)       ;
\end{tikzpicture}
&
\begin{tikzpicture}[baseline,level distance=5mm]
\node(root){\tiny a}[sibling distance=10mm]
        child{node{\tiny b}
                child{node{}
                        child{node{}
                                child{node{}}
                        }
                }
                child{node{}
                        child{node{\tiny c}[sibling distance=6mm]
                                child{node{}}
                                child{node{}}
                        }
                }
        }
        child{node{} edge from parent[draw=none] child{node(right){} edge from
parent[draw=none] child{node(cut0){}}}};
                \fill (cut0) circle (2pt);
                \draw   (root) -- (right)       ;
\end{tikzpicture}

&
\begin{tikzpicture}[baseline,level distance=5mm]
\node(root){\tiny a}[sibling distance=10mm]
        child{node{\tiny b}
                child{node{}
                        child{node{}
                                child{node{}
                                        child{node{}}
                                }
                        }
                }
                child{node{}
                        child{node{\tiny c}[sibling distance=6mm]
                                child{node{\tiny d}
                                        child{node{}}
                                        child{node{}}
                                }
                                child{node(cut1){}}
                        }
                }
        }
        child{node{} edge from parent[draw=none] child{node(right){} edge from
parent[draw=none] child{node(cut0){}}}};
                \fill (cut0) circle (2pt);
                \fill (cut1) circle (2pt);
                \draw   (root) -- (right)       ;
\end{tikzpicture}

&
\begin{tikzpicture}[baseline,level distance=5mm]
\node(root){\tiny a}[sibling distance=10mm]
        child{node{\tiny b}
                child{node{}
                        child{node{}
                                child{node{}
                                        child{node(cut2){}}
                                }
                        }
                }
                child{node{}
                        child{node{\tiny c}[sibling distance=6mm]
                                child{node{\tiny d}
                                        child{node{}
                                                child{node{}}
                                        }
                                        child{node{\tiny e}[sibling
distance=6mm]
                                                child{node{}}
                                                child{node{}}
                                        }
                                }
                                child{node{}}
                        }
                }
        }
        child{node{} edge from parent[draw=none] child{node(right){} edge from
parent[draw=none] child{node(cut0){}}}};
                \fill (cut0) circle (2pt);
                \fill (cut1) circle (2pt);
                \fill (cut2) circle (2pt);
                \draw   (root) -- (right)       ;
\end{tikzpicture}
&
\begin{tikzpicture}[baseline,level distance=5mm]
\node(root){\tiny a}[sibling distance=10mm]
        child{node{\tiny b}
                child{node{}
                        child{node{}
                                child{node{}
                                        child{node{}}
                                }
                        }
                }
                child{node{}
                        child{node{\tiny c}[sibling distance=6mm]
                                child{node{\tiny d}
                                        child{node{}
                                                child{node{}
                                                        child{node{}}
                                                }
                                        }
                                        child{node{\tiny e}[sibling
distance=6mm]
                                                child{node{}
                                                        child{node{}}
                                                }
                                                child{node{}
                                                        child{node{}}
                                                }
                                        }
                                }
                                child{node{}}
                        }
                }
        }
        child{node{} edge from parent[draw=none] child{node(right){} edge from
parent[draw=none] child{node{}}}};
                \fill (cut0) circle (2pt);
                \fill (cut1) circle (2pt);
                \fill (cut2) circle (2pt);
                \draw   (root) -- (right)       ;
\end{tikzpicture}
&
\\\hline
{\tiny\begin{tabular}{c}
\\First\\output tree
\end{tabular}}
&
\begin{tikzpicture}[baseline,level distance=5mm]
\node(root){\tiny a}[sibling distance=6mm]
        child{node{}}
        child{node{}};
\end{tikzpicture}
&
\begin{tikzpicture}[baseline,level distance=5mm]
\node(root){\tiny a}[sibling distance=6mm]
        child{node{}
                child{node{}}
        }
        child{node{} child{node{}}};
\end{tikzpicture}
&
\begin{tikzpicture}[baseline,level distance=5mm]
\node(root){\tiny a}[sibling distance=6mm]
        child{node{}
                child{node{\tiny c}
                        child{node{}}
                        child{node{}}
                }
        }
        child{node{} child{node(cut3){}}};
        \fill (cut3) circle (2pt);
\end{tikzpicture}
&
\begin{tikzpicture}[baseline,level distance=5mm]
\node(root){\tiny a}[sibling distance=6mm]
        child{node{}
                child{node{\tiny c}
                        child{node{\tiny d}
                                child{node{}}
                                child{node{}}
                        }
                        child{node(cut4){}}
                }
        }
        child{node{} child{node{}}};
                \fill (cut3) circle (2pt);
        \fill (cut4) circle (2pt);
\end{tikzpicture}
&
\begin{tikzpicture}[baseline,level distance=5mm]
\node(root){\tiny a}[sibling distance=6mm]
        child{node{}
                child{node{\tiny c}
                        child{node{\tiny d}
                                child{node{}child{node{}}}
                                child{node{}child{node{}}}
                        }
                        child{node(cut4){}}
                }
        }
        child{node{} child{node{}}};
        \fill (cut3) circle (2pt);
        \fill (cut4) circle (2pt);
\end{tikzpicture}
&
\begin{tikzpicture}[baseline,level distance=5mm]
\node(root){\tiny a}[sibling distance=6mm]
        child{node{}
                child{node{\tiny c}
                        child{node{\tiny d}
                                child{node{}child{node{}child{node{}}}}
                                child{node{}child{node{}child{node{}}}}
                        }
                        child{node{}}
                }
        }
        child{node{} child{node{}}};
        \fill (cut3) circle (2pt);
        \fill (cut4) circle (2pt);
\end{tikzpicture}
&
\begin{tikzpicture}[baseline,level distance=5mm]
\node(root){\tiny a}[sibling distance=6mm]
        child{node{} edge from parent[double,solid]
                child{node{\tiny c} edge from parent[double,solid]
                        child{node{\tiny d} edge from parent[double,solid]
                                child{node{} edge from parent[double,solid]
child{node{}  edge from parent[double,solid] child{node{} edge from
parent[double,solid]}}}
                                child{node{}child{node{}child{node{}}}}
                        }
                        child{node{}}
                }
        }
        child{node{} child{node{}}};
        \fill (cut3) circle (2pt);
        \fill (cut4) circle (2pt);
\end{tikzpicture}
\\\hline
{\tiny \begin{tabular}{c}
\\Second\\output tree
\end{tabular}}
&
\begin{tikzpicture}[baseline,level distance=5mm]
\node(root){\tiny b}[sibling distance=6mm]
        child{node{}}
        child{node{}};
\end{tikzpicture}
&
\begin{tikzpicture}[baseline,level distance=5mm]
\node(root){\tiny b}[sibling distance=6mm]
        child{node{}child{node{}}}
        child{node{}child{node{}}};
\end{tikzpicture}
&
\begin{tikzpicture}[baseline,level distance=5mm]
\node(root){\tiny b}[sibling distance=6mm]
        child{node{}child{node{}child{node{}}}}
        child{node{}child{node{}child{node{}}}};
\end{tikzpicture}
&
\begin{tikzpicture}[baseline,level distance=5mm]
\node(root){\tiny b}[sibling distance=6mm]
        child{node{}child{node{}child{node{}child{node{}}}}}
        child{node{}child{node{}child{node{}child{node{}}}}};
\end{tikzpicture}
&
\begin{tikzpicture}[baseline,level distance=5mm]
\node(root){\tiny b}[sibling distance=6mm]
        child{node{}child{node{}child{node{}child{node(cut5){}}}}}
        child{node{}child{node{}child{node{}child{node{\tiny e}
                child{node{}}
                child{node{}}
        }}}};
        \fill (cut5) circle (2pt);
\end{tikzpicture}
&
\begin{tikzpicture}[baseline,level distance=5mm]
\node(root){\tiny b}[sibling distance=6mm]
        child{node{}child{node{}child{node{}child{node{}}}}}
        child{node{}child{node{}child{node{}child{node{\tiny e}
                child{node{}child{node{}}}
                child{node{}child{node{}}}
        }}}};
        \fill (cut5) circle (2pt);
\end{tikzpicture}
&
\begin{tikzpicture}[baseline,level distance=5mm]
\node(root){\tiny b}[sibling distance=6mm]
        child{node{}child{node{}child{node{}child{node{}}}}}
        child{node{} edge from parent[double,solid]child{node{} edge from
parent[double,solid]child{node{} edge from parent[double,solid]child{node{\tiny
e} edge from parent[double,solid]
                child{node{} edge from parent[double,solid]child{node{} edge
from parent[double,solid]}}
                child{node{}child{node{}}}
        }}}};
        \fill (cut5) circle (2pt);
\end{tikzpicture}
\\\hline
Computing & a path &&&&&&\\
{\tiny step 1.} & {\tiny step 2.} & {\tiny step 3.} & {\tiny step 4.} & {\tiny step 5.} & {\tiny step 6.} & {\tiny step 7.}&
\\\hline
\begin{tikzpicture}[baseline,level distance=5mm]
\node(root){\tiny a}[sibling distance=6mm]
        child{node{} edge from parent[double,solid]}
        child{node{}};
\end{tikzpicture}
&
\begin{tikzpicture}[baseline,level distance=5mm]
\node(root){\tiny a}[sibling distance=6mm]
        child{node{} edge from parent[double,solid]}
        child{node{}};
\end{tikzpicture}
&
\begin{tikzpicture}[baseline,level distance=5mm]
\node(root){\tiny a}[sibling distance=6mm]
        child{node{} edge from parent[double,solid]}
        child{node{}};
\end{tikzpicture}
&
\begin{tikzpicture}[baseline,level distance=5mm]
\node(root){\tiny a}[sibling distance=6mm]
        child{node{} edge from parent[double,solid]
                child{node{\tiny c} edge from parent[double,solid]
                        child{node{} edge from parent[double,solid]}
                        child{node{}}
                }
        }
        child{node{} child{node(cut3){}}};
        \fill (cut3) circle (2pt);
\end{tikzpicture}
&
\begin{tikzpicture}[baseline,level distance=5mm]
\node(root){\tiny a}[sibling distance=6mm]
        child{node{} edge from parent[double,solid]
                child{node{\tiny c} edge from parent[double,solid]
                        child{node{\tiny d} edge from parent[double,solid]
                                child{node{} edge from parent[double,solid]}
                                child{node{}}
                        }
                        child{node(cut4){}}
                }
        }
        child{node{} child{node{}}};
                \fill (cut3) circle (2pt);
        \fill (cut4) circle (2pt);
\end{tikzpicture}
&
\begin{tikzpicture}[baseline,level distance=5mm]
\node(root){\tiny a}[sibling distance=6mm]
        child{node{} edge from parent[double,solid]
                child{node{\tiny c} edge from parent[double,solid]
                        child{node{\tiny d} edge from parent[double,solid]
                                child{node{} edge from parent[double,solid]}
                                child{node{}}
                        }
                        child{node(cut4){}}
                }
        }
        child{node{} child{node{}}};
                \fill (cut3) circle (2pt);
        \fill (cut4) circle (2pt);
\end{tikzpicture}
&
\begin{tikzpicture}[baseline,level distance=5mm]
\node(root){\tiny a}[sibling distance=6mm]
        child{node{} edge from parent[double,solid]
                child{node{\tiny c} edge from parent[double,solid]
                        child{node{\tiny d} edge from parent[double,solid]
                                child{node{} edge from parent[double,solid]}
                                child{node{}}
                        }
                        child{node(cut4){}}
                }
        }
        child{node{} child{node{}}};
                \fill (cut3) circle (2pt);
        \fill (cut4) circle (2pt);
\end{tikzpicture}
&
{\tiny\begin{tabular}{c}\\ \\Path\\through the\\first\\output tree\end{tabular}}
\\\hline
\begin{tikzpicture}[baseline,level distance=5mm]
\node(root){\tiny b}[sibling distance=6mm]
        child{node{}}
        child{node{}};
\end{tikzpicture}
&
\begin{tikzpicture}[baseline,level distance=5mm]
\node(root){\tiny b}[sibling distance=6mm]
        child{node{}}
        child{node{} edge from parent[double,solid]};
\end{tikzpicture}
&
\begin{tikzpicture}[baseline,level distance=5mm]
\node(root){\tiny b}[sibling distance=6mm]
        child{node{}}
        child{node{} edge from parent[double,solid]};
\end{tikzpicture}
&
\begin{tikzpicture}[baseline,level distance=5mm]
\node(root){\tiny b}[sibling distance=6mm]
        child{node{}}
        child{node{} edge from parent[double,solid]};
\end{tikzpicture}
&
\begin{tikzpicture}[baseline,level distance=5mm]
\node(root){\tiny b}[sibling distance=6mm]
        child{node{}}
        child{node{} edge from parent[double,solid]};
\end{tikzpicture}
&
\begin{tikzpicture}[baseline,level distance=5mm]
\node(root){\tiny b}[sibling distance=6mm]
        child{node{}}
        child{node{} edge from parent[double,solid]};
\end{tikzpicture}
&
\begin{tikzpicture}[baseline,level distance=5mm]
\node(root){\tiny b}[sibling distance=6mm]
        child{node{}}
        child{node{} edge from parent[double,solid]};
\end{tikzpicture}
&
{\tiny\begin{tabular}{c}Path\\through the\\second\\output tree\end{tabular}}
\\\hline
\begin{tikzpicture}[baseline,level distance=5mm]
\node(root){{\tiny a}}[sibling distance=10mm]
        child{node{{\tiny b}}  edge from parent[double,solid]
                child{node{}}
                child{node{}}
        }
        child{node{} edge from parent[draw=none] child{node(right){} edge from
parent[draw=none]}};
                \draw   (root) -- (right)       ;
\end{tikzpicture}
&
\begin{tikzpicture}[baseline,level distance=5mm]
\node(root){{\tiny a}}[sibling distance=10mm]
        child{node{{\tiny b}}  edge from parent[double,solid]
                child{node{}}
                child{node{} edge from parent[double,solid]}
        }
        child{node{} edge from parent[draw=none] child{node(right){} edge from
parent[draw=none]}};
                \draw   (root) -- (right)       ;
\end{tikzpicture}
&
\begin{tikzpicture}[baseline,level distance=5mm]
\node(root){\tiny a}[sibling distance=10mm]
        child{node{\tiny b}  edge from parent[double,solid]
                child{node{}
                        child{node{}}
                }
                child{node{}  edge from parent[double,solid]
                        child{node{}  edge from parent[double,solid]}
                }
        }
        child{node{} edge from parent[draw=none] child{node(right){} edge from
parent[draw=none] child{node{}}}};
                \draw   (root) -- (right)       ;
\end{tikzpicture}
&
\begin{tikzpicture}[baseline,level distance=5mm]
\node(root){\tiny a}[sibling distance=10mm]
        child{node{\tiny b}  edge from parent[double,solid]
                child{node{}
                        child{node{}
                                child{node{}}
                        }
                }
                child{node{}  edge from parent[double,solid]
                        child{node{\tiny c}[sibling distance=6mm]  edge from
parent[double,solid]
                                child{node{}  edge from parent[double,solid]}
                                child{node{}}
                        }
                }
        }
        child{node{} edge from parent[draw=none] child{node(right){} edge from
parent[draw=none] child{node(cut0){}}}};
                \fill (cut0) circle (2pt);
                \draw   (root) -- (right)       ;
\end{tikzpicture}

&
\begin{tikzpicture}[baseline,level distance=5mm]
\node(root){\tiny a}[sibling distance=10mm]
        child{node{\tiny b}  edge from parent[double,solid]
                child{node{}
                        child{node{}
                                child{node{}
                                        child{node{}}
                                }
                        }
                }
                child{node{}  edge from parent[double,solid]
                        child{node{\tiny c}[sibling distance=6mm]  edge from
parent[double,solid]
                                child{node{\tiny d}  edge from
parent[double,solid]
                                        child{node{}  edge from
parent[double,solid]}
                                        child{node{}}
                                }
                                child{node(cut1){}}
                        }
                }
        }
        child{node{} edge from parent[draw=none] child{node(right){} edge from
parent[draw=none] child{node(cut0){}}}};
                \fill (cut0) circle (2pt);
                \fill (cut1) circle (2pt);
                \draw   (root) -- (right)       ;
\end{tikzpicture}

&
\begin{tikzpicture}[baseline,level distance=5mm]
\node(root){\tiny a}[sibling distance=10mm]
        child{node{\tiny b}  edge from parent[double,solid]
                child{node{}
                        child{node{}
                                child{node{}
                                        child{node(cut2){}}
                                }
                        }
                }
                child{node{}  edge from parent[double,solid]
                        child{node{\tiny c}[sibling distance=6mm]  edge from
parent[double,solid]
                                child{node{\tiny d}  edge from
parent[double,solid]
                                        child{node{}  edge from
parent[double,solid]
                                                child{node{}  edge from
parent[double,solid]}
                                        }
                                        child{node{\tiny e}[sibling
distance=6mm]
                                                child{node{}}
                                                child{node{}}
                                        }
                                }
                                child{node{}}
                        }
                }
        }
        child{node{} edge from parent[draw=none] child{node(right){} edge from
parent[draw=none] child{node(cut0){}}}};
                \fill (cut0) circle (2pt);
                \fill (cut1) circle (2pt);
                \fill (cut2) circle (2pt);
                \draw   (root) -- (right)       ;
\end{tikzpicture}
&
\begin{tikzpicture}[baseline,level distance=5mm]
\node(root){\tiny a}[sibling distance=10mm]
        child{node{\tiny b} edge from parent[double,solid]
                child{node{}
                        child{node{}
                                child{node{}
                                        child{node{}}
                                }
                        }
                }
                child{node{} edge from parent[double,solid]
                        child{node{\tiny c}[sibling distance=6mm] edge from
parent[double,solid]
                                child{node{\tiny d}  edge from
parent[double,solid]
                                        child{node{} edge from
parent[double,solid]
                                                child{node{} edge from
parent[double,solid]
                                                        child{node{} edge from
parent[double,solid]}
                                                }
                                        }
                                        child{node{\tiny e}[sibling
distance=6mm]
                                                child{node{}
                                                        child{node{}}
                                                }
                                                child{node{}
                                                        child{node{}}
                                                }
                                        }
                                }
                                child{node{}}
                        }
                }
        }
        child{node{} edge from parent[draw=none] child{node(right){} edge from
parent[draw=none] child{node{}}}};
                \fill (cut0) circle (2pt);
                \fill (cut1) circle (2pt);
                \fill (cut2) circle (2pt);
                \draw   (root) -- (right)       ;
\end{tikzpicture}
&
{\tiny\begin{tabular}{c}\\ \\Path \\through the\\input tree\end{tabular}}
\\\hline
\end{tabular}
\caption{Illustrating the algorithm underlying Proposition \ref{prop:activevertices}}
\label{fig:exampleactive}
\end{figure}

\end{exa}

As a consequence from the independent choice theorem in \cite{paulybrattka} (or rather its proof) together with $\C_{\Cantor,\sharp=n} \equivW \ic{=n}$ and $\C_{\Cantor,\sharp\leq n} \equivW \ic{\leq n}$ we obtain the following, showing ultimately that picking an element from a finite number of 2-element sets in parallel is just as hard as picking finitely many times from finite sets, with the later questions depending on the answers given so far:
\begin{observation}
\label{obs:hierarchy}
$\ic{=n} \times \ic{=m} \leqW \ic{=n} \star \ic{=m} \leqW \ic{=(nm)}$ and $\ic{\leq n} \times \ic{\leq m} \leqW \ic{\leq n} \star \ic{\leq m} \leqW \ic{\leq (nm)}$
\end{observation}

We can observe that the reductions in Proposition \ref{prop:activevertices} and Observation \ref{obs:hierarchy} are uniform in the natural number parameters. Thus, we can form the corresponding coproducts to obtain the following corollaries.

\begin{cor}
$\ic{=2}^* \equivW \left ( \coprod_{n \in \mathbb{N}} \ic{=n} \right ) \equivW \left ( \coprod_{n, k \in \mathbb{N}} \ic{=n}^{(k)} \right )$
\end{cor}

\begin{cor}
$\ic{\leq 2}^* \equivW \left ( \coprod_{n \in \mathbb{N}} \ic{\leq n} \right ) \equivW \left ( \coprod_{n, k \in \mathbb{N}} \ic{\leq n}^{(k)} \right )$
\end{cor}

\begin{prop}[{\cite{paulyincomputabilitynashequilibria,paulybrattka}}]
$\C_{\{0,1\}}^* \equivW \left ( \coprod_{n \in \mathbb{N}} \C_{\{0, \ldots, n\}} \right ) \equivW \left ( \coprod_{n, k \in \mathbb{N}} \C_{\{0, \ldots, n\}}^{(k)} \right )$
\begin{proof}
The reduction from $\C_{\{0,1\}}^*$ to $\left ( \coprod_{n, k \in \mathbb{N}} \C_{\{0, \ldots, n\}}^{(k)} \right )$ is trivial. By the Independent Choice Theorem \cite[Theorem 7.3]{paulybrattka} we conclude that $\C_{\{0, \ldots, n\}}^{(k)} \leqW \C_{\{0,\ldots,(n+1)^k\}}$, and the proof is uniform in $n$ and $k$. Thus, we may conclude $\left ( \coprod_{n, k \in \mathbb{N}} \C_{\{0, \ldots, n\}}^{(k)} \right ) \leqW \left ( \coprod_{n \in \mathbb{N}} \C_{\{0, \ldots, n\}} \right )$. By \cite[Theorem 32]{paulyincomputabilitynashequilibria} we have $\C_{\{0,\ldots,n\}} \leqW \C_{\{0,1\}}^n$, again with a proof uniform in $n$. This provides the remaining reduction $\left ( \coprod_{n \in \mathbb{N}} \C_{\{0, \ldots, n\}} \right ) \leqW \C_{\{0,1\}}^*$.
\end{proof}
\end{prop}

Whether this property (that sequential uses of some closed choice principle are equivalent to parallel uses) also applies to convex choice $\xc{1}$ remains open at this stage.

\begin{question}
Is there some $k \in \mathbb{N}$ such that $\xc{1} \star \xc{1} \leqW \xc{1}^k$?
\end{question}

The preceding question gains in relevance in light of the following:

\begin{prop}
$\xc{k} \leqW \xc{1}^{(k)}$.
\begin{proof}
In a compact product space, we can compute projections (e.g.~\cite[Proposition 6 (8)]{pauly-synthetic-arxiv}, and projections of convex subsets to the first component are convex themselves. Hence, given the input $A \in \mathcal{A}(\uint^k)$ we can use the first application of $\xc{1}$ to find some $x \in \uint$ such that $\{x\} \times \uint^{k-1} \cap A \neq \emptyset$. But this intersection again is a convex set, and we can use the second application of $\xc{1}$ to find a valid value of the second component, etc. With all $k$ uses of $\xc{1}$, we then obtain a point inside the input set.
\end{proof}
\end{prop}
\subsection{Simplex choice}
The central idea of this subsection is to relate sets of cardinality $n + 1$ to $n$-dimensional convex sets by using the points in the former as the vertices of a simplex. In this, the notion of affine independence features prominently. We remind the reader that points $v_1, \ldots, v_n \in \uint^k$ are called affinely independent, if $\sum_{i=1}^n \lambda_iv_i = 0$ and $\sum_{i=1}^n \lambda_i = 0$ implies $\lambda_j = 0$ for all $1 \leq j \leq n$. Alternative characterizations are that $v_2 - v_1, v_3 - v_1, \ldots, v_n - v_1$ are linearly independent, or that $v_j \notin \operatorname{ConvexHull}(\bigcup_{i=1, i\neq j}^n \{v_j\})$ for all $1 \leq j \leq n$. Note that being affinely independent is an open property, i.e.~making small perturbations to affinely independent points $v_1, \ldots, v_n$ results in points that are affinely independent again.
\begin{prop}
\label{prop:vandermonde}
Given a closed set $A \subseteq \uint$ with $|A| \leq n$, we can compute a closed set $B \subseteq \uint^{n-1}$ with $|A| = |B|$, $\pi_1(B) = A$, and such that the points in $B$ are affinely independent.
\begin{proof}
The function $f:[0,1]\to[0,1]^{n-1}$ defined by $f(x):=(x,x^2,\dots,x^{n-1})$ is computable. As $\uint$ and $\uint^{n-1}$ are computably compact and computably Hausdorff (cf.~\cite{pauly-synthetic-arxiv}), we can compute $f[A] \in \mathcal{A}(\uint^{n-1})$. By definition of $f$, $\pi_1\circ f$ is the identity, in particular $|A| = |f[A]|$ follows.

Let $A \subseteq \{x_1,\dots,x_n\}$ where the $x_i$ are all distinct. Then the determinant of the Vandermonde matrix below is non-zero.

\[\left(\begin{array}{cccc}1 & x_1 & \dots & x_1^{n-1}\\ \vdots &\vdots & &\vdots  \\ 1 & x_n & \dots & x_n^{n-1}\end{array}\right)\]

Subtracting the first row to every other row does not modify the non-zero determinant, but it shows that the vectors $f(x_i)-f(x_1)$ are linearly independent for $2\leq i\leq n$, i.e.~the points in $f[A]$ are indeed affinely independent.
\end{proof}
\end{prop}

\begin{prop}
\label{prop:addpointtosimplex}
Given a closed set $A \subseteq \uint^n$ with $|A| = n + 1$ such that the points in $A$ are affinely independent, we can compute a set $\left (A \cup \{c\}\right )$, where $c$ is a point in the interior of the convex hull of $A$.
\end{prop}

For the proof of this proposition we shall require some preparation. In this, we consider $\mathbb{R}^n$ to be equipped with the Euclidean metric.

\begin{lem}\label{lem:aff-comb}
Let $x_0,x_1,\dots,x_n\in \mathbb{R}^n$ be affinely independent, let $z$ be an affine combination $\sum_{k}\beta_kx_k$ with $\beta_k > 0$ for all $k$, and for all $0 \leq i \leq n$ let $y_i$ be an affine combination $\sum_{k}\alpha_{k,i}x_k$ with $\alpha_{k,i} \leq 0$ for all $k\neq i$. Then $y_0,x_1,\dots,x_n$ are again affinely independent, the affine decomposition $z = \beta'_0 y_0 + \sum_{0 < k}\beta'_k x_k$ satisfies $0 < \beta'_k$ for all $k$, and for all $i > 0$ the affine decomposition $y_i = \gamma_{0,i} y_0 + \sum_{0<k}\gamma_{k,i}x_k$ satisfies $\gamma_{k,i} \leq 0$ for all $k \neq i$.
\begin{proof}
As they are the coefficients of an affine combination, we find that $\sum_k \alpha_{k,i} = 1$ for each $i$, hence $\alpha_{i,i} \geq 1$. The equality $x_0 = \alpha_{0,0}^{-1}y_0 + \sum_{0<k}(-\alpha_{k,0}\alpha_{0,0}^{-1})x_k$ shows that $y_0,x_1,\dots x_n$ are affinely independent. Furthermore $z = \beta_0\alpha_{0,0}^{-1}y_0 + \sum_{0<k}(\beta_k - \beta_0\alpha_{k,0}\alpha_{0,0}^{-1})x_k$, where all the coefficients are positive since $\alpha_{k,0} \leq 0$ for $k \neq 0$, and $y_i = \alpha_{0,i}\alpha_{0,0}^{-1}y_0 + \sum_{0<k}(\alpha_{k,i} - \alpha_{0,i}\alpha_{k,0}\alpha_{0,0}^{-1})x_k$ for all $i$, where $\alpha_{k,i} - \alpha_{0,i}\alpha_{k,0}\alpha_{0,0}^{-1}$ is non-positive for $i \neq 0,k$.
\end{proof}
\end{lem}

By invoking Lemma~\ref{lem:aff-comb} up to $n+1$ times one can prove the following.

\begin{cor}\label{cor:aff-comb}
Let $x_0,x_1,\dots,x_n\in \mathbb{R}^n$ be affinely independent, let $y_i = \sum_{k}\alpha_{k,i}x_i$ with $\alpha_{k,i} \leq 0$ for all $k\neq i$, and let $z$ be an affine combination $\sum_{i}\beta_ix_i$ with $\beta_i > 0$ for all $i$. Then $z$ is an affine combination of the affinely independent $y_0,\dots,y_n$ with positive coefficients.
\end{cor}

\begin{defi}\hfill
\begin{itemize}
\item Given $n + 1$ affinely independent points $x_0, \ldots, x_n\in\mathbb{R}^n$, let $h(x_0, \ldots, x_n)$ be the minimum among the euclidian distances between one point $x_i$ and the hyperplane containing all the other points. We call $h$ the minimal height of $\{x_0, \ldots, x_n\}$.
\item Given $x_0, \ldots, x_n$ be points in a metric space, let $\Delta(x_0,\dots,x_n) : =  \max\{d(x_i,x_j)\,|\, i \neq j\}$.
\item Given $n + 1$ open rational balls $B_0, \ldots, B_n$ in $\uint^n$, we call $$\textrm{JC}(B_0, \ldots, B_n) := \bigcap_{x_0 \in B_0, \ldots, x_n \in B_n} \operatorname{ConvexHull}(\{x_0,\ldots,x_n\})$$ their joint centre.
\end{itemize}
\end{defi}

\begin{lem}\label{lem:aff-comb-shift}
Let $x_0,\dots,x_n\in \mathbb{R}^n$ be affinely independent, let $0 \leq r < (n+1)^{-1}$, let $z = \sum_{0\leq i\leq n}\beta_ix_i$ be an affine combination such that $r < \beta_i$ for all $i$, and let $d(x_i,y_i) < r\cdot h(x_0,\dots,x_n)$ for all $i$. Then $z$ is in the interior of the convex hull of the $y_i$.
\begin{proof}
For all $i$ let $x'_i := (1-nr)x_i + \sum_{k\neq i}rx_k$, such that $x_i = (1-r)(1-(n+1)r)^{-1}x'_i + \sum_{k\neq i}r((n+1)r-1)^{-1} x'_k$ for all $i$, which shows that $x'_0,\dots,x'_n$ are independent. By substitution we find that every affine combination $\sum_{i}\lambda_ix_i$ can be rewritten $\sum_i(\lambda_i-r)(1-(n+1)r)^{-1}x'_i$. For all $i$ let $H_i$ be the hyperplane of the affine combinations $\sum_{k\neq i}\lambda_kx_k$, and $H'_i$ that of the affine combinations $rx_i + \sum_{k\neq i}\lambda_kx_k$. Since $d(H_i,H'_i) = r\cdot d(x_i,H_i) \geq r\cdot h(x_0,\dots,x_n) > d(x_i,y_i)$ and since $\{x_i\} = \cap_{k \neq i} H_k$, for all $i$ the point $y_i = \sum_k\alpha_{k,i} x_i$ satisfies $\alpha_{k,i} < r$ for all $k\neq i$. So $y_i = \sum_k(\alpha_{k,i}-r)(1-(n+1))^{-1}x'_k$ satisfies $\alpha_{k,i}-r < 0$ for all $k\neq i$. Since $z$ is in the interior of the convex hull of the $x'_i$, it is also in that of the $y_i$ by Corollary~\ref{cor:aff-comb}.
\end{proof}
\end{lem}

\begin{lem}\label{lem:aff-comb-shift2}
Let $n \geq 1$ and $x_0,\dots,x_n\in \mathbb{R}^n$ be affinely independent and let some open ball $B$ contain $x_0$ and $x_1$. We can find $\rho > 0$ and $z\in \mathbb{R}^n$ such that $B(z,\frac{\rho}{n}) \subseteq B \cap \textrm{JC}(B(x_0,\rho), \ldots, B(x_n,\rho))$.
\begin{proof}
Let $h := h(x_0,\dots,x_n)$ and $\Delta := \Delta(x_0,\dots,x_n)$. Let $0 < \rho < \frac{h^2}{2(n+1)\Delta}$ be such that the closed balls $\bar{B}(x_0, \frac{3\rho n\Delta}{h})$ and $\bar{B}(x_1,\frac{3\rho n\Delta}{h})$ are included in $B$. Let us consider $L$ the line segment between $x_0$ and the point $\sum_{0<i} n^{-1}x_i$. Let $g$ be the length of $L$ and let $\{z\}$ be its intersection with the sphere $S(x_0,\frac{2\rho n\Delta}{h})$, which is in the interior of $B$ by construction. Note that $h \leq g \leq \Delta$. Let $r := \frac{\rho}{h}$ and let us show that the affine decomposition $z = \sum_{i}\alpha_i x_i$ satisfies $2r \leq \alpha_i$ for all $i$. First, it is easy to check that $\frac{2\rho}{h} \leq 1 - \frac{2\rho n\Delta}{h^2}$, so $2r \leq 1 - \frac{2\rho n\Delta}{hg} = \alpha_0$. Second, $\alpha_i = \alpha_j$ for all $i,j > 0$, since $z \in L$, and $2r = \frac{2\rho}{h} \leq \frac{2\rho n\Delta}{h}\cdot\frac{1}{gn} = \alpha_1$. For all $i$ let $z_i := z + \frac{rd}{n+1}x_i - \sum_{j \neq i}\frac{r}{n+1}x_j$ (so $\{z_0,\dots,z_n\}$ is a rescaling of $\{x_0,\dots,x_n\}$ by factor $r$). The $z_i$ are all in $\bar{B}(x_0, \frac{3\rho n\Delta}{h}) \subseteq B$ since $d(z,z_i) \leq \Delta r \leq \frac{\rho n \Delta}{h}$ for all $i$, and so is their convex hull. Moreover, the affine decomposition of each $z_i$ along the $x_j$ involves coefficients greater than $r$ only, so by Lemma~\ref{lem:aff-comb-shift} the $z_i$ are all in the interior of $\textrm{JC}(B(x_0,\rho), \ldots, B(x_n,\rho))$, and so is the ball $B(z,\frac{\rho}{n})$ inscribed in their convex hull.
\end{proof}
\end{lem}

\begin{proof}[Proof of Proposition \ref{prop:addpointtosimplex}]
  We search for a cover of $A$ by $n + 1$ open rational balls
  $B(x_0,\rho), \ldots,$ $B_n(x_n,\rho)$, with the $x_i$ affinely
  independent and $\rho$ chosen from the $x_i$ as in Lemma
  \ref{lem:aff-comb-shift2}. As $\rho$ depends on the points only via
  their height and mutual distance, any set of $n + 1$ affinely
  independent points is covered by such a collection of
  balls. Moreover, as $\uint^n$ is computably compact, we are
  guaranteed to find a cover eventually.

By Lemma \ref{lem:aff-comb-shift2} we can identify some ball $B(z,t) \subseteq \textrm{JC}(B(x_0,\rho), \ldots, B_n(x_n,\rho))$. Our initial approximation for the output set is $A_0 := \overline{B}(z,t) \cup \bigcup_{i \leq n} \overline{B}(x_i,\rho)$.

We then search for covers of $A$ by smaller and smaller balls contained in the original $B(x_i,\rho)$, i.e.~by some $B(x_i^k,2^{-k}\rho)$ with $B(x_i^k,2^{-k}\rho) \subseteq B(x_i,\rho)$. As the joint center is antimonotone, we find that $B(z,t) \subseteq \textrm{JC}(B(x_0^k,2^{-k}\rho), \ldots, B(x_n^k,2^{-k}\rho))$. We then set our $k$-th approximation to $A_k := \overline{B}(z,2^{-k}t) \cup \bigcup_{i \leq n} \overline{B}(x_i^k,2^{-k}\rho)$.

Now, it may happen finitely many times that in our new cover there are two balls contained in the same $B(x_j,\rho)$. This means that one $B(x_i,\rho)$ is disjoint from $A$ and that the
current $B(z,t)$ could be outside of the convex hull of $A$. In this case, we have to choose some new values for $z', t'$ using Lemma \ref{lem:aff-comb-shift2} again. As we can choose these values such that $B(z',t') \subseteq B(x_j,\rho)$, we find that the resulting approximation is still a subset of the preceding one.

Ultimately, our output is obtained as $\bigcap_{k \in \mathbb{N}} A_k$.
\end{proof}

\begin{prop}
\label{prop:simplexchoice}
Given a finite closed set $A \subseteq \uint^n$, such that the points in $A$ are affinely independent, as well as a point $x$ in the convex hull of $A$, we can compute a point in $A$.
\begin{proof}
Let $p$ be a name of $A=:\{a_0,\dots,a_l\}$ of dimension $l\leq n$ and let $q$ be a name of a point $x$ in the convex hull of $A$. Let us compute one point in $A$ as the intersection of a nested sequence of closed balls $B_k$. We aim at guaranteeing the following conditions for all $k$.
\begin{enumerate}
\item\label{cond:diam} The diameter of $B_k$ is at most $\frac{1}{2^{k-2}}$.
\item\label{cond:cvx-hull} The point $x$ is not in the convex hull of $A\backslash B_k$.
\item\label{cond:bound-inter} The boundary of $B_k$ does not intersect $A$.
\end{enumerate}
Let us define $B_k$ by induction, starting with $B_0:=[-1,2]^n$ so that all conditions are met, and assume that $B_k$ is defined and meets all conditions. Let us enumerate $p$ until there exist $l+1$ disjoint closed balls $b_0,\dots,b_l$ of diameters less than $\frac{1}{2^{k}}$ such that their union includes the approximation of $A$ so far (such balls exist since $|A|=l+1$) and the $b_i$ are fully either inside or outside $B_k$ (thanks to Condition~\ref{cond:bound-inter}).

Since $x$ is in $\operatorname{ConvexHull}(A)$ by assumption, $x=\sum_{i=0}^l\alpha_ia_i$, where $\sum\alpha_i=1$ and $0\leq\alpha_i$. Assume that $x$ is also in the affine span of $A\backslash B_k$, that is, $x=\sum_{i=0}^m\lambda_ia_i$, where $\sum\lambda_i=1$ and, modulo renaming, the $b_{m+1},\dots,b_l$ are exactly the $b_i$ that are included in $B_k$. Since $x=\sum_{i=0}^l\alpha_ia_i=\sum_{i=0}^m\lambda_ia_i$, the $\lambda_i$ are positive by uniqueness of the coefficients of the affine representation, so $x$ is in $\operatorname{ConvexHull}(A\backslash B_k)$, which contradicts Condition~\ref{cond:cvx-hull}. Therefore $x$ is not in the affine span of $A\backslash B_k$, so there exists $j$ such that $x$ is not in the affine span of $A\backslash b_j$. Such a $j$ may be identified in finite time by running in parallel the name $q$ of $x$ and names of the $\operatorname{ConvexHull}(A\backslash b_i)$ that are derived from $p$ (the name provided for $A$). Now let $B_{k+1}$ be a closed ball centered like $b_j$, with a diameter between $\frac{1}{2^{k}}$ and $\frac{1}{2^{k-1}}$, and disjoint from the other $b_i$. Let us show that the three conditions are met. Condition~\ref{cond:diam} is met by construction; Condition~\ref{cond:cvx-hull} is met since $x$ does not belong to $\operatorname{ConvexHull}(A\backslash b_j)$; and Condition~\ref{cond:bound-inter} is also met since $b_j\subseteq B_k$ and $B_{k+1} \cap b_i = \emptyset$ for all $i \neq j$.

The intersection of the $B_k$ is a singleton since their diameters converge to $0$, and each $B_k$ intersects $A$ due to Condition~\ref{cond:cvx-hull}, so $\cap_{k\in\mathbb{N}} B_k=\{a_i\}$ for some $0\leq i\leq l$.
\end{proof}
\end{prop}

\begin{cor}
\label{corr:maincorr}
$\ic{\leq n} \leqW \xc{n-1}$
\begin{proof}
We show $\C_{\uint,\sharp \leq n} \leqW \xc{n-1}$ instead. Given a set of up to $n$ points in $\uint$, we can use Proposition \ref{prop:vandermonde} to turn them into the vertices of a proper simplex in $\uint^{n-1}$. The convex hull is computable as a closed set by Proposition \ref{prop:convexhull}, and we can use $\xc{n-1}$ to pick a point inside the convex hull. Then Proposition \ref{prop:simplexchoice} allows us to recover one of the vertices, which by Proposition \ref{prop:vandermonde} suffices to obtain an element of the original set.
\end{proof}
\end{cor}

\begin{thm}
$\ic{=n} \leW \ic{=n+1}$
\begin{proof}
By Corollary \ref{corr:icnrichcompact}, we can freely change the space we are working in among any rich computably compact computable metric space. We start with an $n$-point subset of $\uint$ and apply Proposition \ref{prop:vandermonde} to obtain $n$ affinely independent points in $\uint^{n-1}$. Then we use Proposition \ref{prop:addpointtosimplex} to obtain a set of cardinality $n+1$ containing the $n$ previous points and some additional point in the interior of their convex hull. This is a valid input to $\ic{=n+1}$ (using Corollary \ref{corr:icnrichcompact} again), and we obtain one of the points, which certainly is contained in the convex hull. Hence, Proposition \ref{prop:simplexchoice} allows us to find one of the vertices, which by Proposition \ref{prop:vandermonde} is sufficient to compute one of the points in the original set.

That the reduction is strict follows from Propositions \ref{prop:cntoicnplusone} and \ref{prop:cntoicn}.
\end{proof}
\end{thm}
Note that while $\ic{\leq n} \leqW \ic{\leq n+1}$ is trivially true, the positive part of the preceding result is not obvious.

\subsection{Application of the large radius technique}
The usefulness of the large radius technique for disproving reducibility to convex choice lies in the observation that convex sets with large outer radius are simpler, as we can then cut by a hyperplane and obtain another convex set of smaller dimension. For convenience, we work with the maximum metric on $\uint^n$ in the following proposition. Choosing another compatible metric may require adapting the precise choice of numbers, yet does not impact the results further along.

\begin{prop}[Cutting]
\label{prop:cutting}
$\xc{n, \lambda > m^{-1}} \leqW \xc{n-1} \star \C_{\{1, \ldots, (m-1)n\}}$
\begin{proof}
We can compute the $(m-1)n$-many hyperplanes $\{im^{-1}\} \times \uint \times \ldots \times \uint$, and $\uint \times \{im^{-1}\} \times \ldots \times \uint$, and so on, and finally $\uint \times \ldots \times \uint \times \{im^{-1}\}$, denoted by $P_{ij}$ with $1 \leq i \leq m - 1$ and $1 \leq j \leq n$. For each of these, we compute the intersection with the given convex set $A$, which will be a convex set itself. As we are in a compact space, we can detect emptiness, in particular, we can compute $\{(i, j) \mid P_{ij} \cap A \neq \emptyset\} \in \mathcal{A}(\{1,\ldots, n\} \times \{1, \ldots, m-1\}) \cong \mathcal{A}(\{1, \ldots, (m-1)n\})$. The guarantee $\lambda(A) > m^{-1}$ implies that for some $\langle i, j\rangle$ we have $P_{ij} \cap A \neq \emptyset$. Application of $\C_{\{1, \ldots, (m-1)n\}}$ allows us to find a suitable pair $(i,j)$. Then we compute the projection (which is possible, again, as we are in a compact space, \cite[Proposition 6(8)]{pauly-synthetic-arxiv}) of $P_{ij} \cap A$ to the components distinct from $j$, which will be a non-empty convex subset of a $n - 1$ dimensional space, and given a point from the latter convex set, by inserting $im^{-1}$ as the $j$-th component, we obtain a point in $A$.
\end{proof}
\end{prop}

\begin{cor}
\label{corr:largediameterpluscutting}
Let $\C_{\mathbf{X}}|_{\mathfrak{A}}$ and $(\C_{\mathbf{X}})_{X_\varepsilon(\mathfrak{A})}$ be fractals\footnote{Cf.~Open question \ref{question:fractals}.} and $\C_{\mathbf{X}}|_{\mathfrak{A}} \leqW \xc{n+1}$. Then $(\C_{\mathbf{X}})_{X_\varepsilon(\mathfrak{A})} \leqW \xc{n}$ for all $\varepsilon > 0$.
\begin{proof}
Corollary \ref{corr:ldpfractals} gives us $(\C_{\mathbf{X}})_{X_\varepsilon(\mathfrak{A})} \leqW \xc{n+1, \lambda > m^{-1}}$ for some $m \in \mathbb{N}$, then Proposition \ref{prop:cutting} implies $(\C_{\mathbf{X}})_{X_\varepsilon(\mathfrak{A})} \leqW \xc{n} \star \C_{\{1, \ldots, (m-1)(n+1)\}}$, and finally Theorem \ref{theo:fractalabsorption} fills the gap to $(\C_{\mathbf{X}})_{X_\varepsilon(\mathfrak{A})} \leqW \xc{n}$.
\end{proof}
\end{cor}

For $n \geq k \geq 1$ let $\ic{=n \rhd k} := \C_{\uint }|_{ \{A \in \mathcal{A}(\uint) \mid |A| = n \wedge |\{i < n \mid  [\frac{2i}{2n}, \frac{2i+1}{2n}] \cap A \neq \emptyset\}| \geq k\}}$. So $\ic{=n \rhd k}$ is choice for $n$ element sets, where we know that our set intersects at least $k$ of a collection of fixed distinct regions. We shall need three properties of these choice principles:

\begin{prop}
\label{prop:icboxes}
\begin{enumerate}
\item $\ic{=n \rhd (k+1)} \leqW (\C_{\uint})_{X_{(5n)^{-1}}(\dom(\ic{=n \rhd k}))}$.
\item $\ic{=n+1 \rhd n}$ is not computable.
\item Any $\ic{=n \rhd k}$ is a fractal.
\end{enumerate}
\begin{proof}
\begin{enumerate}
\item We use a representation $\psi_-$ of $\dom(\ic{=n \rhd (k+1)})$ where the finite approximation available at any stage lies in the interior of the one available at the previous stage (\cite[Proposition 3.4]{paulybrattka3}). Any such name already belongs to  $X_{(5n)^{-1}}(\dom(\ic{=n \rhd k}))$. To see this, note that any ball of radius $(5n)^{-1}$ can intersect at most one of the $k + 1$ inhabited regions for $\ic{=n \rhd (k+1)}$, hence removing such a ball leaves at least $k$ regions inhabited. It remains to split some of the remaining approximations of points into several to keep the cardinality condition satisfied, but this is unproblematic.
\item Any algorithm solving $\ic{=n+1 \rhd n}$ would need to eventually pick on of the regions. However, when we represent the sets with names where the approximation at any finite stage is in the interior of the approximation at the previous stage (as obtained by \cite[Proposition 3.4]{paulybrattka3}), we can then make sure that the selected region contains two points. But then the algorithm would have to solve $\ic{=2}$, contradiction.
\item We use the same representation $\psi_-$ of $\dom(\ic{=n \rhd k})$ as in (1). Restricting $\ic{=n \rhd k} \circ \psi_-$ to some arbitrary clopen set is equivalent to picking some arbitrary finite approximation $A$ and restricting $\ic{=n \rhd k}$ to subsets of $A$. For any corresponding finite approximation $A$, we consider those regions $[\frac{2i}{2n}, \frac{2i+1}{2n}]$ where $\left ([\frac{2i}{2n}, \frac{2i+1}{2n}] \cap A\right)^\circ \neq \emptyset$. Let $I$ be the set of corresponding indices. For $i \in I$ we then pick a non-degenerate rational interval $[a_i, b_i] \subseteq [\frac{2i}{2n}, \frac{2i+1}{2n}] \cap A$. Let $B = \bigcup_{i \in I} [a_i, b_i]$.

    We shall use $i^+ := \min \left ( \{j \in I \mid j > i\} \cup \{n\} \right )$ and $i^- := \max \left (\{j \in I \mid j < i\} \cup \{-1\} \right )$. Moreover, we understand $a_n = 1$.

    Now we shall argue that $\ic{=n \rhd k} \leqW \ic{=n \rhd k}|_{\{A \in \mathcal{A}(\uint) \mid A \subseteq B\}}$, which establishes $\ic{=n \rhd k}$ as a fractal. As any $p \in \dom(\ic{=n \rhd k} \circ \psi_-)$ encodes a finite set, we will for each $i \in I$ eventually detect some rational $c_i \notin \psi_-(p)$ with $b_i \leq \frac{2i+1}{2n} < c_i < \frac{2i^+}{2n} \leq a_{i^+}$. Let $c_{0} = 0$. For each $i \in I$, we can now pick a computable homeomorphism $R_i : [c_{i^-}, c_i] \to [a_i, b_i]$. Next, we join all $R_i$ to obtain a computably invertible computable map $R : \uint \setminus \{c_i \mid i \in I\} \to B$ and proceed to use $R$ to rescale $\psi_-(p)$ to a subset of $B$. By construction we see that $\{i \mid \psi_-(p) \cap [\frac{2i}{2n}, \frac{2i+1}{2n}] \neq \emptyset\} \subseteq \{i \mid R[\psi_-(p)] \cap [\frac{2i}{2n}, \frac{2i+1}{2n}] \neq \emptyset\}$, hence the procedure does indeed map inputs for $\ic{=n \rhd k}$ to inputs for $\ic{=n \rhd k}|_{\{A \in \mathcal{A}(\uint) \mid A \subseteq B\}}$. As $R$ is computably invertible, we can recover a valid output to the former from any valid output of the latter. Thus, the reduction is demonstrated.
\end{enumerate}
\end{proof}
\end{prop}

\begin{cor}
\label{corr:step}
$\ic{=n\rhd k} \leqW \xc{m+1}$ implies $\ic{=n \rhd (k+1)} \leqW \xc{m}$.
\begin{proof}
$\ic{=n\rhd k} \leqW \xc{m+1}$ implies $(\C_{\uint})_{X_{(5n)^{-1}}(\dom(\ic{=n \rhd k}))} \leqW  \xc{m+1}|_{\lambda > l^{-1}}$ for suitable $l \in \mathbb{N}$ by Proposition \ref{prop:icboxes} (3) and Corollary \ref{corr:ldpfractals}. Using Proposition \ref{prop:icboxes} (1) on the left-hand side, and Proposition \ref{prop:cutting} on the right-hand side yields $\ic{=n \rhd (k+1)} \leqW \xc{m} \star \C_{\{1, \ldots, (m+1)l\}}$. Then Proposition \ref{prop:icboxes} (3) together with Theorem \ref{theo:fractalabsorption} provides the claim.
\end{proof}
\end{cor}

\begin{thm}
\label{main:theo}
$\ic{=n+2} \nleqW \xc{n}$.
\begin{proof}
Assume $\ic{=n+2} \equivW \ic{=(n+2) \rhd 1} \leqW \xc{n}$. Iterated use of Corollary \ref{corr:step} allows us to conclude that $\ic{=(n+2) \rhd (n+1)} \leqW \xc{0}$, i.e.~$\ic{=(n+2) \rhd (n+1)}$ is computable, which contradicts Proposition \ref{prop:icboxes} (2).
\end{proof}
\end{thm}

\begin{cor}
\label{corr:convexhierarchy}
$\xc{n} \leW \xc{n+1}$.
\begin{proof}
Combine Theorem \ref{main:theo} with Corollary \ref{corr:maincorr}.
\end{proof}
\end{cor}
\subsection{Beyond compact spaces}
\label{subsec:beyondcompactness}
In this subsection we shall investigate $\C_{\mathbb{R}^k,\sharp \leq n}$, $\C_{\mathbb{R}^k,\sharp = n}$ and $\C_{\mathbb{R}^k }|_{ \{A \in \mathcal{A}(\mathbb{R}^k) \mid A \textnormal{ is convex}\}}$. Essentially, all these choice principles are in the same relation to the corresponding ones for compact spaces as the full $\C_\mathbb{R}$ has to $\C_\uint$. We recall from \cite{paulybrattka} that $\C_\mathbb{R} \equivW \C_\uint \star \C_\mathbb{N}$ and $\C_{\mathbb{R},\sharp = 1} \equivW \C_\mathbb{N}$.

\begin{prop}
\label{prop:ricn}
For any $k, n \geq 1$ it holds that $\C_{\mathbb{R}^k,\sharp = n} \equivW \C_\mathbb{N}$.
\begin{proof}
The proof of $\C_{\mathbb{R}^k,\sharp = n} \leqW \C_\mathbb{N}$ is similar to the proof of Proposition \ref{prop:icnleqwcn}: We guess $n$ disjoint rational hypercubes encoded as some $m \in \mathbb{N}$, and reject the guess if the intersection of any of these hypercubes (as compact sets) and the input set is empty. If the rational hypercubes are a suitable guess, then all points are available as compact singletons, hence can be computed.

For the other direction, we make use of $\C_\mathbb{N} \equivW \C_{\mathbb{N},\sharp=1}$ from \cite{paulybrattka}. Given some closed singleton $A \in \mathcal{A}(\mathbb{N})$, we can compute the set $\{(i + \frac{j}{n+1}, 0, \ldots, 0) \mid i \in A \wedge 1 \leq j \leq n\} \in \mathcal{A}(\mathbb{R}^k)$, and any point from the latter set suffices to reconstruct $i \in \mathbb{N}$.
\end{proof}
\end{prop}

\begin{cor}
$\C_{\mathbb{R}^k,\sharp = n} \equivW  \ic{=n} \star \C_\mathbb{N}$.
\begin{proof}
By the independent choice theorem (\cite{paulybrattka}), we have $\C_{\mathbb{N}} \star \C_{\mathbb{N}} \equivW \C_{\mathbb{N}}$. With Proposition \ref{prop:icnleqwcn}, we then see $\ic{=n} \star \C_\mathbb{N} \equivW \C_\mathbb{N}$.
\end{proof}
\end{cor}

\begin{prop}
$\C_{\mathbb{R}^k,\sharp \leq n} \equivW \ic{\leq n} \star \C_\mathbb{N}$
\label{prop:icleqnr}
\begin{proof}
Given some $A \in \mathcal{A}(\mathbb{R}^k)$, we can compute \[\{\langle d_1, \ldots, d_k\rangle \mid A \cap \left ( [d_1, d_1 +1] \times \ldots \times [d_k, d_k+1] \right ) \neq \emptyset\}\,,\] hence we can use $\C_\mathbb{N}$ to find some compact hypercube having non-empty intersection with the input set. This intersection clearly satisfies the cardinality restriction, thus is a suitable input for $\ic{\leq n}$, showing $\C_{\mathbb{R}^k,\sharp \leq n} \leqW \ic{\leq n} \star \C_\mathbb{N}$.

For the other direction we assume w.l.o.g.~that $k = 1$. As in the proof of Proposition \ref{prop:ricn}, we use $\C_\mathbb{N} \equivW \C_{\mathbb{N},\sharp=1}$ from \cite{paulybrattka}. Moreover, we use $\C_{\uint,\sharp \leq n}$ in place of $\ic{\leq n}$ (Corollary \ref{corr:icnrichcompact}). In the call to $\C_{\uint,\sharp \leq n} \star \C_{\mathbb{N},\sharp=1}$, let $A$ be the input used for $\C_{\mathbb{N},\sharp=1}$, and let $A_i$ be the set used as input to $\C_{\uint,\sharp \leq n}$ if $i$ is obtained as answer from $\C_{\mathbb{N},\sharp=1}$, or the empty set, if $i$ is not a valid output. Now it is possible to compute $\left ( \bigcup_{i \in A} (\{2i\} + A_i)\right ) \in \mathcal{A}(\mathbb{R})$, which is a suitable input to $\C_{\mathbb{R},\sharp \leq n}$ and any element of this set encodes all needed information to solve the instance to $\C_{\uint,\sharp \leq n} \star \C_{\mathbb{N},\sharp=1}$.
\end{proof}
\end{prop}

\begin{prop}
$\C_{\mathbb{R}^k }|_{ \{A \in \mathcal{A}(\mathbb{R}^k) \mid A \textnormal{ is convex}\}} \equivW \xc{k} \star \C_\mathbb{N}$
\begin{proof}
The same reduction witness that is used in the proof of Proposition \ref{prop:icleqnr} to show $\C_{\mathbb{R}^k,\sharp \leq n} \leqW \ic{\leq n} \star \C_\mathbb{N}$ works for $\C_{\mathbb{R}^k }|_{ \{A \in \mathcal{A}(\mathbb{R}^k) \mid A \textnormal{ is convex}\}} \leqW \xc{k} \star \C_\mathbb{N}$, too, as the intersection of a convex set and a hypercube is a convex set. The other direction, again, proceeds exactly as in Proposition \ref{prop:icleqnr}.
\end{proof}
\end{prop}

The joint structure of the preceding propositions generally demonstrates that for any class of closed subsets $\mathfrak{A}$ that is closed under rescaling and either closed under intersection with intervals, or only contains bounded sets, it follows that $\C_{\mathbb{R}}|_{\mathfrak{A}} \equivW \C_{\uint}|_{\mathfrak{A}} \star \C_\mathbb{N}$. Such results may fruitfully interplay with Theorem \ref{theo:closedfractalabsorption}.

\section{Finding zeros of functions with finitely many local extrema}
As a closed subset of a computable metric space can equivalently be expressed as the zero set of some continuous function into $\mathbb{R}$, we recognize $\ic{=n}$ ($\ic{\leq n}$) to simultaneously be the degree of finding a zero of a function on a rich computably compact computable metric space, in particular a function $f : \uint \to \mathbb{R}$ that has exactly (up to) $n$ zeros. However, usually when such a task is encountered, the bound on the number of zeros is linked to a bound on the number of local extrema (we understand this to exclude plateaus and the end points of the interval).

We shall now demonstrate that the restriction to a bounded number of local extrema makes the search for zeros significantly easier. The underlying algorithmic result is that given such a function, we can compute a fixed finite number of real numbers that will include all zeros of the function at hand.

\begin{thm}
\label{theo:boundedminima}
For any $n \in \mathbb{N}$ the multivalued map $\operatorname{Zero}_{n \min} : \subseteq \mathcal{C}([a,b], \mathbb{R}) \mto [a,b]^{3^n}$ is computable, where $f \in \dom(\operatorname{Zero}_{n \min})$ iff $f$ has up to $n$ local minima and $f(a) \neq 0 \neq f(b)$, and $(x_1, \ldots, x_{3^n}) \in \operatorname{Zero}_{n \min}(f)$ iff $f(x) = 0$ implies $\exists i . x_i = x$.
\begin{proof}
Our algorithm proceeds in a divide-and-conquer method, subdividing the interval into smaller intervals such that the corresponding restrictions of the function have fewer local extrema. At each stage, we have some number of guesses for potential zeros available. We always start with the assumption that our function has at most one zero in the current interval, and give approximations for all guesses accordingly. If this assumption fails, it has to do so through the detectable existence of certain obstructions. The identification of an obstruction allows us to produce two restrictions of the function, and to allocate our guesses to the restrictions in a valid manner.

We keep track of a bound on the number of local minima at each stage. More precisely, we only consider \emph{essential} extrema, which are those that could be linked to zeros. If $f(a) > 0$, then the left-most local maximum is inessential (this may be at $a$ itself), likewise $f(a) < 0$ makes the left-most infimum inessential, $f(b) > 0$ the right-most maximum and $f(b) < 0$ the right-most minimum.

By switching to $-f$ instead if necessary, we can restrict ourselves to the situation where $f(a) > 0$, and distinguish the cases $f(b) > 0$ and $f(b) < 0$. Note that in both settings, the number of essential maxima cannot exceed the number of essential minima. Hence, we can retain the bound for the number of essential minima when switching the sign.

We describe the number of guesses required (given some bound on the number of essential minima) by two functions $\kappa_1, \kappa_2 : \mathbb{N} \to \mathbb{N}$, one for each of the two situation. These functions are described via a recurrence relation, which in turn is found by investigating by how much the subdivision reduces the bound on the number of essential minima. By the following Lemma \ref{lemma:recurrence}, we see that we need more potential solutions in the second situation, and that the number of potential solutions our algorithm may need to provide is given by $3^n$.

{\bf Situation 1}: $f(b) > 0$, $j$ essential minima

\begin{figure}[htbp]
\centering
\includegraphics[width=0.5\textwidth]{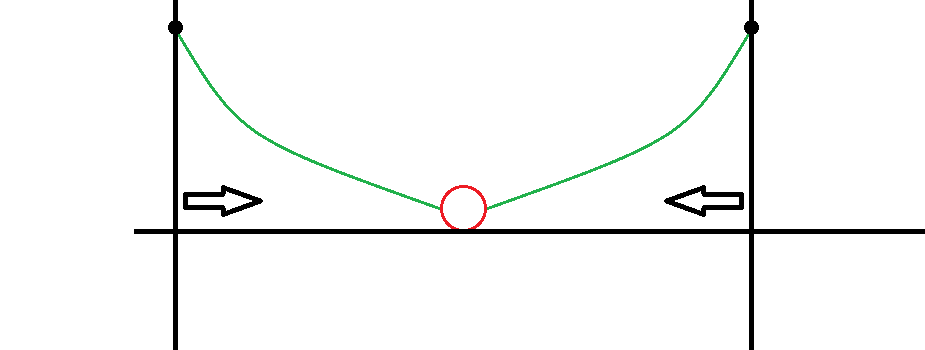}
\caption{The base case in the situation $f(b) > 0$}
\label{base1}
\end{figure}

We start by shrinking the interval from both sides, as long as we can establish that the function is strictly positive in the area. Furthermore, we pay attention to the configurations discussed below (also Figure \ref{obs11}, \ref{obs12}), which may block the shrinkage. If any such configuration exists, it will eventually be found, and if they are all absent, the process will collapse the interval to a singleton, hence assigning all available solution attempts to this value. In this case, the function takes only strictly positive values outside of the remaining point, hence, this point is the only potential zero of the function. If the function has no (essential) minima, it cannot have a zero, hence we may set $\kappa_1(0) := 0$.

\begin{figure}[htbp]
\centering
\includegraphics[width=0.5\textwidth]{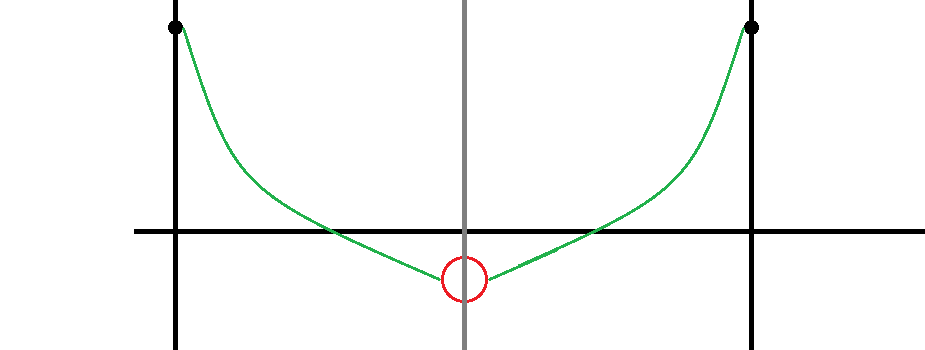}
\caption{The first obstruction in the situation $f(b) > 0$}
\label{obs11}
\end{figure}

If we find some $x \in [a,b]$ such that $f(x) < 0$ can be proven, we split the interval into the parts $[a, x]$ and $[x, b]$. This split renders one minimum inessential, but we do not know in which part the remaining essential minima and maxima will end up. Both parts belong to Situation 2 (described below), yielding the (partial) recurrence relation $\kappa_1(j+1) \geq 2\kappa_2(j)$ that we need to enforce.

\begin{figure}[htbp]
\centering
\includegraphics[width=0.5\textwidth]{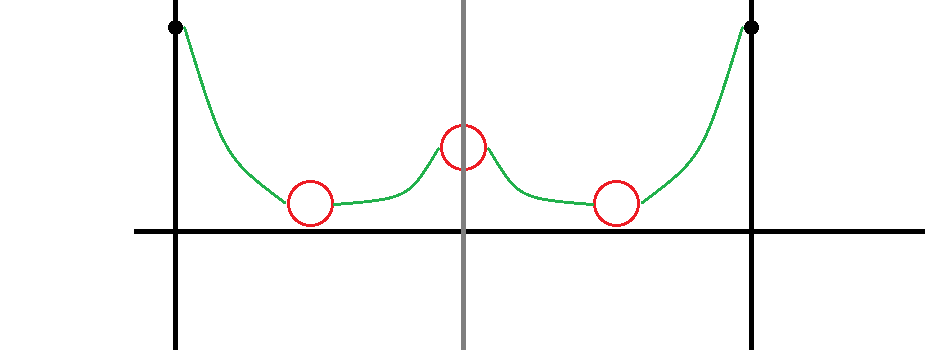}
\caption{The second obstruction in the situation $f(b) > 0$}
\label{obs12}
\end{figure}

If we find $x < y < z$ such that $f(x) < f(y) > f(z)$ and $f(y) > 0$ can be proven, we split into $[a, y]$ and $[y, b]$. Here a local maximum is rendered inessential, and furthermore we know that there is at least one essential local minimum in each part, both parts belong to the first situation. Hence the second part of the recurrence relation is $\kappa_1(j+1) \geq 2\kappa_1(j)$. Together with the first inequality, we conclude that setting $\kappa_1(j+1) := \max \{2\kappa_1(j), 2\kappa_2(j)\}$ fulfills the requirements.

{\bf Situation 2}: $f(b) < 0$, $j$ essential minima

\begin{figure}[htbp]
\centering
\includegraphics[width=0.5\textwidth]{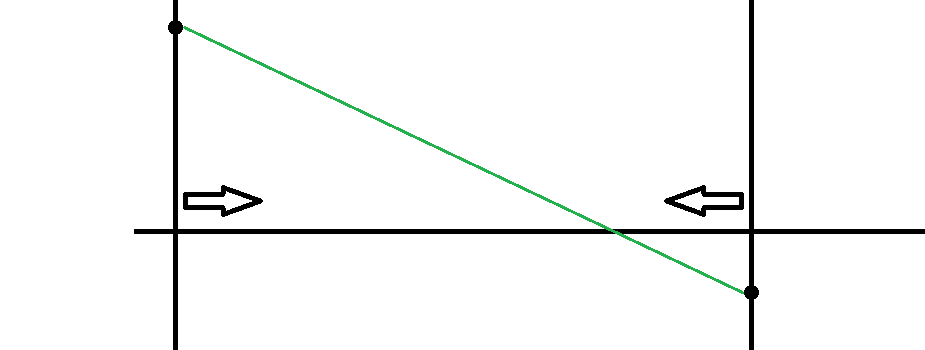}
\caption{The base case in the situation $f(b) < 0$}
\label{base2}
\end{figure}

Again the interval is shrunk from both sides, on the left as long as the function is known to be strictly positive, on the right as long as the function is known to be strictly negative. If none of the configurations listed below (also Figures \ref{obs21}, \ref{obs22}) is ever detected, the interval will collapse to a single point, which is the unique zero of the function in the interval. As both obstructions require $j > 0$, we can set $\kappa_2(0) := 1$.

\begin{figure}[htbp]
\centering
\includegraphics[width=0.5\textwidth]{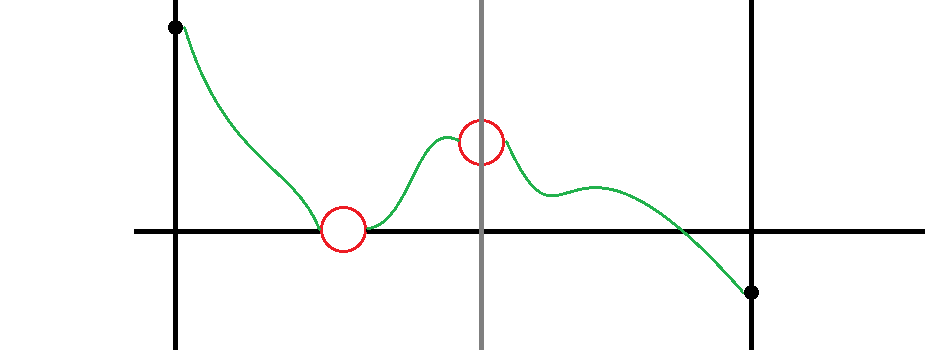}
\caption{The first obstruction in the situation $f(b) < 0$}
\label{obs21}
\end{figure}

If we find $x < y$ such that $f(x) < f(y) > 0$ can be proven, we split into $[a, y]$ and $[y, b]$. A maximum is rendered inessential, and we know that at least one essential minimum is in the left part. The left part belongs to the first situation, and the right part to the second. The corresponding inequality is $\kappa_2(j+1) \geq \kappa_2(j) + \kappa_1(j+1)$.

\begin{figure}[htbp]
\centering
\includegraphics[width=0.5\textwidth]{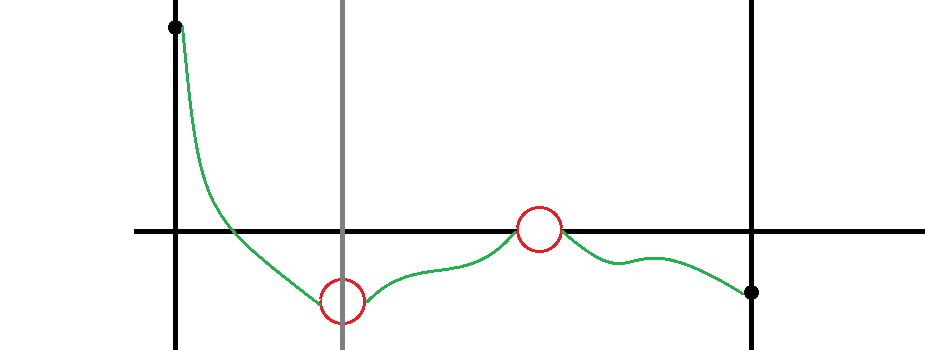}
\caption{The second obstruction in the situation $f(b) < 0$}
\label{obs22}
\end{figure}

If we find $x < y$ such that $0 > f(x) < f(y)$ can be proven, we split into $[a, x]$ and $[x, b]$. A minimum is rendered inessential. The left part belongs to the second situation, and the right part to the first (after moving to $-f$). The corresponding term is again $\kappa_2(j+1) \geq \kappa_1(j+1) + \kappa_2(j)$, hence our recurrence relation is $\kappa_2(j+1) := \kappa_1(j+1) + \kappa_2(j)$.
\end{proof}
\end{thm}

\begin{lem}
\label{lemma:recurrence}
The following recurrence relation \[\begin{array}{rcl} \kappa_1(0) & = & 0 \\ \kappa_2(0) & = & 1 \\ \kappa_1(j+1) & = & \max \{2\kappa_1(j), 2\kappa_2(j)\} \\ \kappa_2(j+1) & = & \kappa_1(j+1) + \kappa_2(j)\end{array}\] has the solution $\kappa_1(j+1) = 2*3^j$ and $\kappa_2(j) = 3^j$.
\begin{proof}
From the last equation we deduce $\kappa_2(j) \geq \kappa_1(j)$; thus the third equation simplifies to $\kappa_1(j+1) = 2\kappa_2(j)$. This in turn renders the last equation into $\kappa_2(j+1) = 3\kappa_2(j)$. Together with the second equation, we then may conclude $\kappa_2(j) = 3^j$. The claim for $\kappa_1$ now follows immediately from the (simplified) third equation.
\end{proof}
\end{lem}

\begin{cor}
\label{corr:boundedminima}
Finding a zero of a continuous function $f : \uint \to \mathbb{R}$ with $f(0) > 0 \wedge f(1) \neq 0$ and up to $n$ local minima is Weihrauch reducible to $\C_{\{1, \ldots, 3^n\}}$.
\begin{proof}
Using Theorem \ref{theo:boundedminima} we compute the $3^n$-many potential solutions, and then test for each of them whether the function value actually is zero at that point. If the function value at the $i$-th solution is recognized to be non-zero, $i$ is removed from the input to $\C_{\{1, \ldots, 3^n\}}$. Any number remaining is the index of a zero of the function.
\end{proof}
\end{cor}

\begin{cor}
Finding a zero of a continuous function $f : \uint \to \mathbb{R}$ with $f(0) > 0 \wedge f(1) \neq 0$ and up to $m > 1$ zeros is not reducible to finding a zero of a function $g : \uint \to \mathbb{R}$ with $g(0) > 0 \wedge g(1) \neq 0$ and up to $n$ local minima for any $n > 0$.
\begin{proof}
Combine Corollary \ref{corr:boundedminima} and Corollary \ref{corr:icnnleqc1n}.
\end{proof}
\end{cor}

In the remainder of this section we shall demonstrate that the algorithm in Theorem \ref{theo:boundedminima} is optimal in two senses: First, it is not possible to reliably compute $3^n-1$ real numbers containing all zeros of a suitable function. Second, using closed choice for a finite space to actually find a zero as in Corollary \ref{corr:boundedminima} is optimal in the sense of Weihrauch reducibility (up to the precise cardinality of the space involved).

\begin{thm}
For no $n \in \mathbb{N}$ the multivalued map $\operatorname{Zero}'_{n \min} : \subseteq \mathcal{C}(\uint, [-1,1]) \mto \uint^{3^n-1}$ is computable, where $f \in \dom(\operatorname{Zero}'_{n \min})$ iff $f$ has up to $n$ local minima and $f(0) > 0 > f(1)$, and $(x_1, \ldots, x_{3^n-1}) \in \operatorname{Zero}'_{n \min}(f)$ iff $f(x) = 0$ implies $\exists i . x_i = x$.
\begin{proof}
Assume the contrary, and let $N$ be the least $n$ for which the map $\operatorname{Zero}'_{n \min}$ is computable. We can exclude $N = 0$, as a function without local minima can still have a zero (see Figure \ref{base2}), so $\operatorname{Zero}'_{0 \min}$ is not even well-defined.

Now we will describe an input $f$ for  $\operatorname{Zero}'_{N \min}$ designed to fool the hypothetical algorithm. As $f$ is a continuous function between compact Hausdorff spaces, we can equivalently represent $f$ via its graph as a closed or compact subset of $\uint \times [-1,1]$. This in turn we can assume to be given by coverings of the graph by finitely many rational closed boxes of smaller and smaller total area, such that the projections to the $x$-axis of the boxes involved in a single approximation intersect at at most one point.

\begin{figure}[htbp]
\centering
\includegraphics[width=0.5\textwidth]{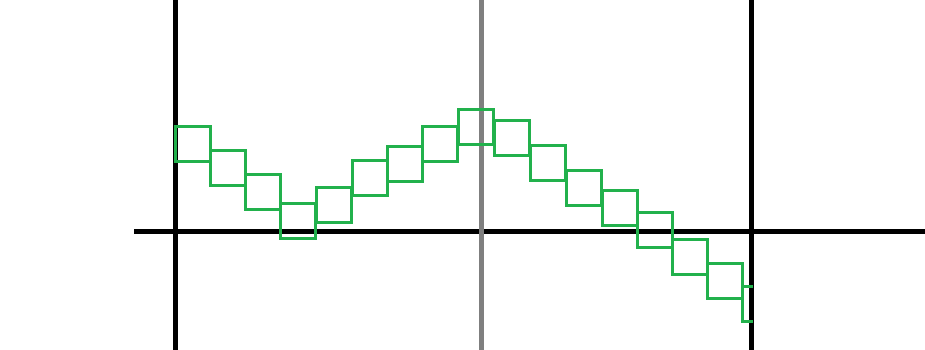}
\caption{The construction of a counterexample}
\label{ce1}
\end{figure}

We start by providing approximations following Figure \ref{ce1}, i.e.~we determine that $f(0.5) > 0$, that there is a local minima to the left of $0.5$ while keeping a box present containing some $[a,b] \times \{0\}$ for $a < b < 0.5$, and while not providing any additional information about the placement of local minima. Note that all zeros of functions admitting a name extending such an approximation are included in one of two disjoint intervals, one to the left and one to the right of $0.5$. In particular, the hypothetical algorithm will have to decide eventually for each of the $3^N - 1$ numbers it is producing to which of these intervals (if any at all) it is going to belong.

Let us assume that less than $3^{N-1}$ numbers are assigned to the right interval. Then for any continuous function $g : \uint \to [-1,1]$ with up to $N-1$ local minima and $g(0) > 0 > g(1)$, we could use rescaling to fit it into the box intersecting the $x$-axis to the right of $0.5$, and then piecewise linear functions to extend it to a function $f$ compatible with the current approximation that has up to $N$ local minima. By providing those numbers as potential solutions that are assigned to the right of the interval, we would have demonstrated that $\operatorname{Zero}'_{N-1 \min}$ is computable, which contradicts our choice of $N$ as the least such number. Hence, at least $3^{N-1}$ numbers have to be assigned to the right interval, leaving at most $2*3^{N-1}-1$ for the left interval.

\begin{figure}[htbp]
\centering
\includegraphics[width=0.5\textwidth]{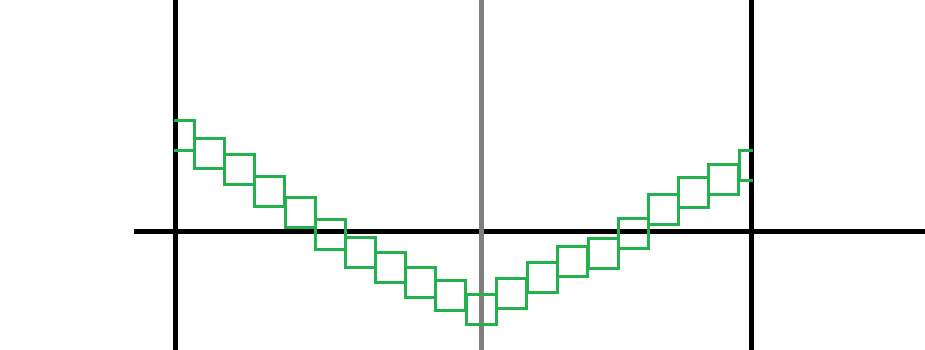}
\caption{The auxiliary construction of a counterexample}
\label{ce2}
\end{figure}

Now we can use the box intersecting the $x$-axis to the left of $0.5$ to provide approximations according to Figure \ref{ce2}, while using piecewise linear functions to extend this to some function $f$ having all its local minima inside the box shown in Figure \ref{ce2}. Again, we see that all zeros in Figure \ref{ce2} are contained in one of two boxes, and our hypothetical algorithm has to decide for each of its $2*3^{N-1}-1$ remaining outputs whether to assign it to the left or the right box. If it assigns less that $3^{N-1}$ numbers to the left box, we can, as before, rescale a function $g : \uint \to [-1,1]$ with up to $N-1$ local minima and $g(0) > 0 > g(1)$ into the right box, and use an otherwise piecewise linear extension adding a single further local minima to conclude computability of $\operatorname{Zero}'_{N-1 \min}$.

Thus, we see that only $3^{N-1}-1$ numbers can be assigned to the right box. But then we could again take a function $g : \uint \to [-1,1]$ with up to $N-1$ local minima and $g(0) > 0 > g(1)$, move to $h(x) := g(1-x)$, rescale, extend in a piecewise linear way and get computability of $\operatorname{Zero}'_{N-1 \min}$. Thus, we cannot avoid the contradiction involved in the assumption $\operatorname{Zero}'_{N \min}$ were computable.
\end{proof}
\end{thm}

\begin{prop}
\label{prop:cnpolynomials}
$\C_{\{1, \ldots, n\}}$ is Weihrauch reducible to finding roots of polynomials of degree $2n$.
\begin{proof}
Note that there is a computable multivalued function $r : \mathbb{S} \mto \mathbb{R}$ such that $r(\top) \subseteq \{x \in \mathbb{R} \mid x > 0\}$ and $r(\bot) = \{0\}$. Given some closed set $A \in \mathcal{A}(\{1, \ldots, n\})$ we can compute the polynomial $p_A := \prod_{i \leq n} \left ((x - \frac{i}{n})^2 + r(A(i))\right )$. Now $\frac{i}{n}$ is a root of $p_A$ iff $i \in A$.
\end{proof}
\end{prop}

\begin{cor}
The following are Weihrauch equivalent:
\begin{enumerate}
\item $\coprod_{n \in \mathbb{N}} \C_{\{1, \ldots, n\}}$
\item Finding a root of a polynomial of known degree $n > 1$
\item Finding a zero of a continuous function $f : \uint \to \mathbb{R}$ with $f(0) > 0 \wedge f(1) \neq 0$ with a known upper bound $n > 0$ on the number of local minima.
\end{enumerate}
\end{cor}

\section{Some related work}
The degree of $\ic{\leq 2}$ was amongst the first Weihrauch degrees to be studied in some detail: In \cite{weihrauchb}, Weihrauch had shown that this map is equivalent to the multivalued function mapping real numbers to an expansion in base $b > 1$. Essentially, the non-computability here lies solely in the fact that some rational numbers have two expansions, whereas the other real numbers have one.

The principle $\ic{\leq 2}$ is also connected to the study on computability on the space of bottomed sequences, or Plotkin's $\mathbb{T}^\mathbb{N}$, pioneered by Tsuiki \cite{tsuiki,tsuiki2,tsuiki3}. In this space, digits may be remain undetermined ($\bot$) for a while -- and potentially for ever -- or be specified as either $0$ or $1$.

\begin{defi}
Let the represented space $\mathbb{T}$ have the underlying set $\{0, 1, \bot\}$ and the representation $\delta_\mathbb{T} : \Baire \to \{0, 1, \bot\}$ be defined by $\delta_\mathbb{T}(0^\mathbb{N}) = \bot$ and $\delta_\mathbb{T}(p) = \min \{n \in \mathbb{N} \mid p(n) \neq 0\} \mod 2$ for $p \neq 0^\mathbb{N}$. Let $\mathbb{T}^\mathbb{N}_m$ be the subspace of $\mathbb{T}^\mathbb{N}$ where at most $m$ components take the value $\bot$.
\end{defi}

These spaces characterize the dimension of computable Polish spaces. Whether there is a connection to the characterization of the dimension of $\uint^k$ via the Weihrauch degree of connected choice; or to the results on dimension and the properties of representations in \cite{pauly-kihara-arxiv}, seems to be an open question.

\begin{thm}[Ohta, Tsuiki, Yamada \cite{tsuiki2}]
A computable Polish spaces $\mathbf{X}$ embeds into $\mathbb{T}^\mathbb{N}_m$ iff $\dim(\mathbf{X}) \leq m$.
\end{thm}

We shall now consider the multivalued functions $\textrm{Concretize}_m : \mathbb{T}_m^\mathbb{N} \mto \Cantor$ defined via $p \in \textrm{Concretize}(q)$ iff $\forall n \in \mathbb{N} \ (q(n) \neq \bot) \Rightarrow (q(n) = p(n))$. At CCC 2014, Brattka raised the question what the Weihrauch degree of $\textrm{Concretize}_m$ would be.

\begin{prop}
$\left (\ic{\leq 2} \right )^m \equivW \textrm{Concretize}_m \equivW
\left (\textrm{Concretize}_1 \right )^m$
\end{prop}\newpage
\begin{proof}\hfill
\begin{description}
\item[$\textrm{Concretize}_m \leqW \left (\ic{\leq 2} \right )^m$]

The proof is using the same idea as the one of Proposition \ref{prop:activevertices}. We may consider the input of $\textrm{Concretize}_m$ to consist of a sequence $(w_i)_{i \in \mathbb{N}}$ of finite words over $\{0,1,\bot\}$ of increasing length such that whenever $w_i(j) \in \{0,1\}$, then $w_i(j) = w_k(j)$ for all $k \geq i$, and moreover, each $w_i$ contains the symbol $\bot$ exactly $m$ times. The output must agree with each of the $w_i$ on the location of the $0$s and $1$s.

Each occurrence of $\bot$ in $w_1$ is associated with some binary tree with two vertices at each layer below the root. As long as subsequent $w_i$ share the same occurrence of $\bot$, both vertices on the current level receive a child in the corresponding tree. If a $\bot$ is overwritten by a $0$ (by a $1$), the right (the left) branch dies out, and the left (the right) branch splits. A new $\bot$ will appear at the end of the word, and will be associated with the same tree as the overwritten one was.

Knowing an infinite branch through each of the $m$ trees then allows us to replace each $\bot$ in the $w_i$ by either $0$ or $1$, depending on the direction the path takes at the corresponding branching vertex. The result are longer and longer prefixes of some valid output.

\item[$\left (\ic{\leq 2} \right )^m \leqW \left (\textrm{Concretize}_1 \right )^m$]

Clearly it suffices to show $\ic{\leq 2} \leqW \textrm{Concretize}_1$. For this, we note that given some binary tree $T$ we can compute a sequence $p \in \mathbb{T}^\mathbb{N}$ where $p(n) = 0$ if the right subtree of the $n$-th vertex dies out first, $p(n) = 1$ if the left subtree of the $n$-th vertex dies out first, and $p(n) = \bot$ if both the left and the right subtree of the $n$-th vertex are infinite. If $T$ is in $\dom(\ic{\leq 2})$, then $p \in \mathbb{T}_1^\mathbb{N}$. Moreover, any $q \in \textrm{Concretize}_1(p)$ computes an infinite path through $T$.

\item[$\left (\textrm{Concretize}_1 \right )^m \leqW \textrm{Concretize}_m$]

The reduction is obtained by some standard bijection between $\left (\mathbb{T}^\mathbb{N}\right )^m$ and $\mathbb{T}^\mathbb{N}$ on the input, and a matching one between $\left (\Cantor\right )^m$ and $\Cantor$ on the output.\qedhere
\end{description}
\end{proof}

\noindent Building upon the conference version of the present paper (\cite{paulyleroux-cie}), Neumann \cite{eike-neumann} has classified the Weihrauch degree of the Browder-G\"ohde-Kirk fixed point theorem. One of his results is that finding a fixed point of a 1-Lipschitz function on $\uint^k$ is Weihrauch equivalent to $\xc{k}$. Together with the results from \cite{paulybrattka3}, this shows a trichotomy for the difficulty of finding fixed points depending on the Lipschitz constant: For $L < 1$, it is computable, for $L = 1$, it is convex choice, and for $L > 1$, it is connected choice. We also point out that as a consequence of Corollary \ref{corr:convexhierarchy}, the Browder-G\"ohde-Kirk fixed point theorem becomes strictly harder with increasing dimension; unlike Brouwer's fixed point theorem, which eventually stabilizes.

\bibliographystyle{eptcs}
\bibliography{../../spieltheorie}

\section*{Acknowledgements}
This work benefited from the Royal Society International Exchange Grant IE111233 and the Marie Curie International Research Staff Exchange Scheme \emph{Computable
Analysis}, PIRSES-GA-2011- 294962. The first author was also supported by the the ERC inVEST (279499) project.

We would like to thank Vladik Kreinovich for a question asked at CCA 2011 that spawned our investigation of the choice principles $\ic{\leq n}$ and $\ic{= n}$. We are indebted to Vasco Brattka for the formation of this line of research, and for providing the one dimensional cases of Propositions \ref{prop:addpointtosimplex}, \ref{prop:simplexchoice}. Proposition \ref{prop:cnpolynomials} was inspired by a remark of Martin Ziegler. The paper benefited from the careful attention of the anonymous referees.

\end{document}